\documentclass{article}%
\usepackage{amsmath}
\usepackage{amsfonts}
\usepackage{amssymb}
\usepackage{graphicx}%
\date{May 8, 2003}
\newtheorem{theorem}{Theorem}[section]
\newtheorem{corollary}[theorem]{Corollary}
\newtheorem{lemma}[theorem]{Lemma}

\newtheorem{remark}[theorem]{Remark}

\newenvironment{proof}[1][Proof]{\noindent\textbf{#1.} }{\ \rule{0.5em}{0.5em}
\medskip\smallskip
}
\numberwithin{equation}{section}
\begin{document}

\title{Higher Order Quasiconvexity Reduces to Quasiconvexity}
\author{Gianni Dal Maso
\and Irene Fonseca
\and Giovanni Leoni
\and Massimiliano Morini}
\maketitle

\begin{abstract}
In this paper it is shown that higher order quasiconvex functions suitable in
the variational treatment of problems involving second derivatives may be
extended to the space of all matrices as classical quasiconvex functions.
Precisely, it is proved that a smooth strictly $2$-quasiconvex function with
$p$-growth at infinity, $p>1$, is the restriction to symmetric matrices of a
$1$-quasiconvex function with the same growth. As a consequence, lower
semicontinuity results for second-order variational problems are deduced as
corollaries of well-known first order theorems.

\end{abstract}

\section{Introduction}

In recent years there has been a renewed interest in higher order variational
problems motivated by various mathematical models in engineering and materials
science: in connection with the so-called gradient theories of phase
transitions within elasticity regimes (see \cite{CFL}, \cite{KM}, \cite{Mue});
in the study of equilibria of micromagnetic materials where mastery of second
order energies (here accounting for the exchange energy) is required (see
\cite{CKO}, \cite{DS}, \cite{Mue}, \cite{RS}); in the theory of second order
structured deformations (SOSD) (see \cite{OP}), in the Blake-Zisserman model
for image segmentation in computer vision (see \cite{CLT}); etc..

In the study of lower semicontinuity, relaxation and $\Gamma$%
-convergence\ problems for second order functional the natural notion of
convexity, $2$-quasiconvexity, was introduced by Meyers in \cite{M} (see also
\cite{BCO}, \cite{F}). We recall that a real valued function $f$, defined on
the space $\mathbb{M}_{\operatorname*{sym}}^{n\times n}$ of ${n\times n}$
symmetric matrices, is \textit{$2$-quasiconvex} if
\[
\int_{Q}\left[  f\left(  A+\nabla^{2}\phi\right)  -f\left(  A\right)  \right]
\,dx\geq0
\]
for every $A\in\mathbb{M}_{\operatorname*{sym}}^{n\times n}$ and every
$\phi\in C_{c}^{2}\left(  Q\right)  $, where $Q:=\left(  0,1\right)  ^{n}$ is
the unit cube.

While lower semicontinuity properties of functionals depending \textit{only}
on second order derivatives can be proved easily, when lower order terms are
present, the question is significantly more difficult, since sufficient tools
to handle localization and truncation of gradients are still missing.

To bypass these difficulties one would be tempted to transform higher order
into first order problems, where one uses the standard notion of
quasiconvexity, called $1$-quasiconvexity in this paper. We recall that a real
valued function $f$, defined on the space $\mathbb{M}^{n\times n}$ of
${n\times n}$ matrices, is \textit{$1$-quasiconvex} if
\[
\int_{Q}\left[  f\left(  A+\nabla\varphi\right)  -f\left(  A\right)  \right]
\,dx\geq0
\]
for every $A\in\mathbb{M}^{n\times n}$ and every $\varphi\in C_{c}^{1}\left(
Q;\mathbb{R}^{n}\right)  $.

Thus we are led to the following question.

\begin{enumerate}
\item[\textbf{(Q)} ] \textit{Is every }$2$\textit{-quasiconvex function the
restriction of a }$1$\textit{-quasiconvex function to the space of symmetric
matrices?}
\end{enumerate}

A good indication of the plausibility of an affirmative answer is that it
holds for polyconvex functions as noticed by Dacorogna and Fonseca. Indeed if
\[
f\left(  A\right)  =g\left(  M\left(  A\right)  \right)  \qquad A\in
\mathbb{M}_{\operatorname*{sym}}^{n\times n},
\]
where $g$ is a convex function and $M\left(  A\right)  $ stands for the vector
whose components are all the minors of $A$, then the function
\[
F\left(  A\right)  :=g\left(  \frac{M\left(  A\right)  +M\left(  A\right)
^{t}}{2}\right)  \qquad A\in\mathbb{M}^{n\times n}%
\]
is a polyconvex extension of $f$ to the whole space $\mathbb{M}^{n\times n}$
of ${n\times n}$ matrices.

It is known that $2$-gradient Young measures, i.e. Young measures generated by
second order gradients, may be characterized by duality via Jensen's
inequality with respect to $2$-quasiconvex functions (see \cite{FM}), just as
gradient Young measures are characterized by duality with $1$-quasiconvex
functions (see \cite{KP}). Therefore, the understanding of the structure of
$2$-gradient Young measures helps the study of $2$-quasiconvex functions, and,
accordingly, the following result by \v{S}ver\'{a}k in \cite[Lemma 1]{S}
provides further evidence that $1$-quasiconvexity and $2$-quasiconvexity are
somehow strictly linked: If a Young measure $\nu$ on $\mathbb{M}^{n\times n}$
is generated by a sequence $\left\{  \nabla u_{k}\right\}  $ of gradients,
with $\left\{  u_{k}\right\}  $ bounded in $W^{1,p}\left(  \Omega
;\mathbb{R}^{n}\right)  $ for some $p>1$, and $\operatorname*{supp}\nu
_{x}\subset\mathbb{M}_{\operatorname*{sym}}^{n\times n}$ for $\mathcal{L}^{n}$
a.e. $x\in\Omega$, then $\nu$ is generated also by a sequence $\left\{
\nabla^{2}w_{k}\right\}  $, with $\left\{  w_{k}\right\}  $ bounded in
$W^{2,p}\left(  \Omega\right)  $.

A partial answer to \textbf{(Q)} was given by M\"{u}ller and \v{S}ver\'{a}k
(see \cite{MS}). Indeed, as an auxiliary result to construct a counter-example
to regularity for elliptic systems, they proved that any smooth, strictly
$2$-quasiconvex function $f:\mathbb{M}_{\operatorname*{sym}}^{2\times
2}\rightarrow\mathbb{R}$, with bounded second derivatives, is the restriction
of a $1$-quasiconvex function. The main purpose of this paper is to extend
their result to any space dimension and to a larger class of strictly
$2$-quasiconvex functions with $p$-growth at infinity, with $p>1$.

\begin{theorem}
\label{theorem 2}Let $f\in C^{1}\left(  \mathbb{M}_{\operatorname*{sym}%
}^{n\times n}\right)  $ satisfy the following conditions for suitable
constants $p>1$, $\mu\geq0$, $L\geq\nu>0$:

\begin{enumerate}
\item[(a) ] (strict $2$-quasiconvexity)%
\[
\int_{Q}\left[  f\left(  A+\nabla^{2}\phi\right)  -f\left(  A\right)  \right]
\,dx\geq\nu\int_{Q}\left(  \mu^{2}+\left\vert A\right\vert ^{2}+\left\vert
\nabla^{2}\phi\right\vert ^{2}\right)  ^{\!\!\frac{p-2}{2}}\left\vert
\nabla^{2}\phi\right\vert ^{2}\,dx
\]
for every $A\in\mathbb{M}_{\operatorname*{sym}}^{n\times n}$ and every
$\phi\in C_{c}^{2}\left(  Q\right)  $;

\item[(b)] (Lipschitz condition for gradients)%
\begin{equation}
\left\vert \nabla f\left(  A+B\right)  -\nabla f\left(  A\right)  \right\vert
\leq L\left(  \mu^{2}+\left\vert A\right\vert ^{2}+\left\vert B\right\vert
^{2}\right)  ^{\frac{p-2}{2}}\left\vert B\right\vert \label{3}%
\end{equation}
for every $A$, $B\in\mathbb{M}_{\operatorname*{sym}}^{n\times n}$.
\end{enumerate}

\noindent Then there exists a $1$-quasiconvex function $F:\mathbb{M}^{n\times
n}\rightarrow\mathbb{R}$ such that
\begin{align}
&  F\left(  A\right)  =f\left(  A\right)  \qquad\forall A\in\mathbb{M}%
_{\operatorname*{sym}}^{n\times n},\label{5}\\
&  \left\vert F\left(  A\right)  \right\vert \leq c_{f}\left(  1+\left\vert
A\right\vert ^{p}\right)  \qquad\forall A\in\mathbb{M}^{n\times n}, \label{6}%
\end{align}
for a suitable constant $c_{f}>0$ depending on $f$.
\end{theorem}

We remark that a $1$-quasiconvex function $F$ satisfying (\ref{5}) and
(\ref{6}) is constructed explicitly if $p\geq2$ (see (\ref{900})), while in
the case $1<p<2$ it is defined as the quasiconvex envelope of a suitable
extension of $f$ to $\mathbb{M}^{n\times n}$.

The proof of the theorem relies on a Korn-type inequality for divergence-free
vector fields and uses heavily the Lipschitz condition on the gradient of $f$.
The use of Korn-type inequalities prevents us from obtaining a similar result
for the case $p=1$, which, if valid, will require a different treatment.

We do not know if the result continues to hold without assuming (\ref{3}).
However, when condition (\ref{3})\ is dropped we can still prove the following
weaker version of Theorem \ref{theorem 2}.

\begin{theorem}
\label{theorem 1}Let $f:\mathbb{M}_{\operatorname*{sym}}^{n\times
n}\rightarrow\mathbb{R}$ satisfy the following conditions for suitable
constants $p>1$, $\mu\geq0$, $\nu>0$, $M>0$:

\begin{enumerate}
\item[(a) ] (strict $2$-quasiconvexity)%
\[
\int_{Q}\left[  f\left(  A+\nabla^{2}\phi\right)  -f\left(  A\right)  \right]
\,dx\geq\nu\int_{Q}\left(  \mu^{2}+\left\vert A\right\vert ^{2}+\left\vert
\nabla^{2}\phi\right\vert ^{2}\right)  ^{\!\!\frac{p-2}{2}}\left\vert
\nabla^{2}\phi\right\vert ^{2}\,dx
\]
for every $A\in\mathbb{M}_{\operatorname*{sym}}^{n\times n}$ and every
$\phi\in C_{c}^{2}\left(  Q\right)  $;

\item[(b)] (growth condition)%
\begin{equation}
\left\vert f\left(  A\right)  \right\vert \leq M\left(  1+\left\vert
A\right\vert ^{p}\right)  \label{8}%
\end{equation}
for every $A\in\mathbb{M}_{\operatorname*{sym}}^{n\times n}$.
\end{enumerate}

\noindent Then there exists an increasing sequence $\left\{  F_{k}\right\}  $
of $1$-quasiconvex functions $F_{k}:\mathbb{M}^{n\times n}\rightarrow
\mathbb{R}$ such that
\begin{align}
&  \lim_{k\rightarrow\infty}F_{k}\left(  A\right)  =f\left(  A\right)
\qquad\forall A\in\mathbb{M}_{\operatorname*{sym}}^{n\times n},\label{1}\\
&  \left\vert F_{k}\left(  A\right)  \right\vert \leq c_{k}\left(
1+\left\vert A\right\vert ^{p}\right)  \qquad\forall A\in\mathbb{M}^{n\times
n}, \label{7}%
\end{align}
for a suitable sequence of constants $\left\{  c_{k}\right\}  $ depending only
on $k$ and on the structural constants $p$, $\mu$, $\nu$, $M$, but not on the
specific function $f$.
\end{theorem}

Theorem \ref{theorem 1} allows us to reduce lower semicontinuity problems for
$2$-quasi\allowbreak convex normal integrands of the form $f=f(x,u,\nabla
u,\nabla^{2}u)$ to first order problems (see Section \ref{section lower} for
more details). Indeed, as a consequence of Theorem \ref{theorem 1} we can
prove the following result, which extends to the second order setting a lower
semicontinuity property of $1$-quasiconvex functions in $SBV(\Omega
;{\mathbb{R}}^{d})$ due to Ambrosio \cite{AFP} and later generalized by
Kristensen \cite{K}. For the definition and properties of the space
$SBH(\Omega)$ we refer to \cite{CLT0} and~\cite{CLT}.

\begin{theorem}
\label{theorem4}Let $\Omega\subset{\mathbb{R}}^{n}$ be a bounded open set and
let
\[
f:\Omega\times{\mathbb{R}}\times{\mathbb{R}}^{n}\times\mathbb{M}%
_{\operatorname*{sym}}^{n\times n}\rightarrow\lbrack0,+\infty)
\]
be an integrand which satisfies the following conditions:

\begin{enumerate}
\item[(a)] the function $f(x,\cdot,\cdot,\cdot)$ is lower semicontinuous on
${\mathbb{R}}\times{\mathbb{R}}^{n}\times\mathbb{M}_{\operatorname*{sym}%
}^{n\times n}$ for $\mathcal{L}^{n}$ a.e. $x\in\Omega$;

\item[(b)] the function $f(x,u,\xi,\cdot)$ is $2$-quasiconvex on
$\mathbb{M}_{\operatorname*{sym}}^{n\times n}$ for $\mathcal{L}^{n}$ a.e.
$x\in\Omega$ and every $\left(  u,\xi\right)  \in{\mathbb{R}}\times
{\mathbb{R}}^{n}$;

\item[(c)] there exist a locally bounded function $a:\Omega\times{\mathbb{R}%
}\times{\mathbb{R}}^{n}\rightarrow\lbrack0,+\infty)$ and a constant $p>1$ such
that%
\[
0\leq f(x,u,\xi,A)\ \leq a\left(  x,u,\xi\right)  (1+|A|^{p})
\]
for $\mathcal{L}^{N}$ a.e.\ $x\in\Omega$ and every $(u,\xi,A)\in{\mathbb{R}%
}^{n}\times{\mathbb{R}}^{n}\times\mathbb{M}_{\operatorname*{sym}}^{n\times n}$.
\end{enumerate}

\noindent Then
\[
\int_{\Omega}f(x,u,\nabla u,\nabla^{2}u)\,dx\leq\liminf_{j\rightarrow\infty
}\int_{\Omega}f(x,u_{j},\nabla u_{j},\nabla^{2}u_{j})\,dx
\]
for every $u\in SBH(\Omega)$ and any sequence $\{u_{j}\}\subset SBH(\Omega)$
converging to $u$ in $W^{1,1}(\Omega)$ and such that
\begin{equation}
\sup_{j}\left(  \left\Vert \nabla^{2}u_{j}\right\Vert _{L^{p}}+\int_{S(\nabla
u_{j})}\theta(\left\vert \left[  \nabla u_{j}\right]  \right\vert
)\,d\mathcal{H}^{n-1}\right)  <\infty, \label{700}%
\end{equation}
where $\theta:\left[  0,\infty\right)  \rightarrow\left[  0,\infty\right)  $
is a concave, nondecreasing function such that%
\[
\lim_{t\rightarrow0^{+}}\frac{\theta\left(  t\right)  }{t}=\infty,
\]
and $\left[  \nabla u_{j}\right]  $ denotes the jump of $\nabla u_{j}$ on the
jump set $S(\nabla u_{j})$.
\end{theorem}

An analogous result has been proved in \cite{FLP} in the space $BH(\Omega)$ in
the case where (\ref{700}) is replaced by
\[
\sup_{j}\left\Vert \nabla^{2}u_{j}\right\Vert _{L^{p}}<\infty\qquad\left\vert
D_{s}^{2}u_{j}\right\vert (\Omega)\rightarrow0,
\]
where $D_{s}^{2}u_{j}$ is the singular part of the $\mathbb{M}%
_{\operatorname*{sym}}^{n\times n}$-valued measure $D^{2}u_{j}$. Note that
condition (\ref{700}) arises naturally in the context of free-discontinuity
problems (see~\cite{AFP}).

As a corollary of Theorem \ref{theorem4} we have the following result.

\begin{corollary}
\label{corollary lower}Let $\Omega$ and $f$ be as in Theorem \ref{theorem4}.
Then
\[
\int_{\Omega}f(x,u,\nabla u,\nabla^{2}u)\,dx\leq\liminf_{j\rightarrow\infty
}\int_{\Omega}f(x,u_{j},\nabla u_{j},\nabla^{2}u_{j})\,dx
\]
for every $u\in W^{2,p}(\Omega)$ and any sequence $\{u_{j}\}\subset
W^{2,p}(\Omega)$ weakly converging to $u$ in $W^{2,p}(\Omega)$.
\end{corollary}

In this generality Corollary \ref{corollary lower} was proved in \cite{FLP}
and under stronger hypotheses in \cite{F}, \cite{GP}, and \cite{M}.

\begin{remark}
\emph{All of the above are still valid in a vectorial setting, i.e.,}
$u:\Omega\rightarrow\mathbb{R}^{d}$, \emph{with} $\mathbb{M}%
_{\operatorname*{sym}}^{n\times n}$ \emph{replaced now by} $\left(
\mathbb{M}_{\operatorname*{sym}}^{n\times n}\right)  ^{d}$. \emph{The proofs
are entirely similar to those presented in this paper for the case} $d=1$,
\emph{and we leave the obvious adaptations to the reader. }
\end{remark}

The paper is organized as follows. In Section \ref{section auxiliary} we
present some auxiliary results including the Korn-type inequality mentioned
above. In Section \ref{section proofs} we prove Theorems \ref{theorem 2} and
\ref{theorem 1}, while Theorem \ref{theorem4} is addressed in the last section.

\section{Auxiliary results\label{section auxiliary}}

We begin with some results on the Helmholtz Decomposition and on Korn's type inequalities.

A function $w:\mathbb{R}^{n}\rightarrow\mathbb{R}^{d}$ is said to be
\textit{$Q$--periodic }if $w(x+e_{i})=w(x)$ for a.e. $x\in\mathbb{R}^{n}$ and
every $i=1,\ldots,n$, where $(e_{1},\ldots,e_{n})$ is the canonical basis of
$\mathbb{R}^{n}$. The spaces of $Q$-periodic functions of
$W_{\operatorname*{loc}}^{1,p}\left(  \mathbb{R}^{n};\mathbb{R}^{n}\right)  $,
$W_{\operatorname*{loc}}^{2,p}\left(  \mathbb{R}^{n}\right)  $, and
$C^{\infty}(\mathbb{R}^{n};\mathbb{R}^{n})$ are denoted by
$W_{\operatorname*{per}}^{1,p}\left(  Q;\mathbb{R}^{n}\right)  $,
$W_{\operatorname*{per}}^{2,p}\left(  Q\right)  $, and $C_{\operatorname*{per}%
}^{\infty}(Q;\mathbb{R}^{n})$, respectively.

\begin{lemma}
[Helmholtz decomposition]\label{lemma 6}For every $p>1$ and every $\varphi\in
W_{\operatorname*{per}}^{1,p}\left(  Q;\mathbb{R}^{n}\right)  $ there exist
two functions $\phi\in W_{\operatorname*{per}}^{2,p}\left(  Q\right)  $ and
$\psi\in W_{\operatorname*{per}}^{1,p}\left(  Q;\mathbb{R}^{n}\right)  $ such
that
\[
\varphi=\nabla\phi+\psi,\quad\operatorname*{div}\psi=0.
\]
The function $\psi$ is uniquely determined, while $\phi$ is determined up to
an additive constant.
\end{lemma}

\begin{proof}
Since, by periodicity, $\operatorname*{div}\varphi$ has zero average on $Q$,
there exists a $Q$-periodic\ solution $\phi$ of the equation $\Delta
\phi=\operatorname*{div}\varphi$, which is unique up to an additive constant.
It is clear now that $\psi:=\varphi-\nabla\phi$ is $Q$-periodic and
$\operatorname*{div}\psi=0$.
\end{proof}

Throughout the paper, for every $A\in\mathbb{M}^{n\times n}$, we denote the
symmetric and antisymmetric parts of $A$ by
\[
A^{s}:=\frac{A+A^{t}}{2},\qquad A^{a}:=\frac{A-A^{t}}{2},
\]
where $A^{t}$ is the transpose matrix of $A$.

If $\psi:\Omega\subset\mathbb{R}^{n}\rightarrow\mathbb{R}^{d}$ is any
function, then $\nabla\psi$\ is a $d\times n$ matrix in $\mathbb{M}^{d\times
n}$, with $\left(  \nabla\psi\right)  _{ij}:=\frac{\partial\psi_{i}}{\partial
x_{j}}.$\ Also, differential operators applied to matrix-valued fields are
understood on a row-by-row basis, e.g., if $\psi:\Omega\subset\mathbb{R}%
^{n}\rightarrow\mathbb{R}^{n}$, then
\[
\operatorname*{div}\nabla\psi:=\left(
\begin{array}
[c]{c}%
\operatorname*{div}\nabla\psi_{1}\\
\vdots\\
\operatorname*{div}\nabla\psi_{n}%
\end{array}
\right)  .
\]
To simplify the notation, for any $\psi:\Omega\subset\mathbb{R}^{n}%
\rightarrow\mathbb{R}^{n}$, we set%
\[
\nabla\psi^{t}:=\left(  \nabla\psi\right)  ^{t},\qquad\nabla\psi^{s}:=\left(
\nabla\psi\right)  ^{s},\qquad\nabla\psi^{a}:=\left(  \nabla\psi\right)  ^{a}.
\]

\begin{lemma}
\label{lemma 1}For every $p>1$ there exists a constant $\gamma_{n,p}\geq1$
such that
\begin{equation}
\int_{Q}\left\vert \nabla\psi\right\vert ^{p}\,dx\leq\gamma_{n,p}\int
_{Q}\left\vert \nabla\psi^{a}\right\vert ^{p}\,dx \label{451}%
\end{equation}
for every $Q$-periodic function $\psi:\mathbb{R}^{n}\rightarrow\mathbb{R}^{n}$
of class $C^{\infty}$ with $\operatorname*{div}\psi=0$.
\end{lemma}

\begin{proof}
Since $\operatorname*{div}\nabla\psi^{t}=\nabla\left(  \operatorname*{div}%
\psi\right)  =0$ we have%
\begin{equation}
\Delta\psi=2\operatorname*{div}\left(  \frac{\nabla\psi-\nabla\psi^{t}}%
{2}\right)  =2\operatorname*{div}\left(  \nabla\psi^{a}\right)  . \label{901}%
\end{equation}
Hence (\ref{451}) follows from standard $L^{p}$ estimates for periodic
solutions of the Poisson equation (see \cite{GT}).
\end{proof}

Next we study the behavior of auxiliary functions of the type%
\begin{equation}
g\left(  x\right)  :=\left(  \mu^{2}+\left\vert x\right\vert ^{2}\right)
^{\frac{p}{2}} \label{450}%
\end{equation}
defined on an arbitrary Hilbert space $X$.

\begin{lemma}
\label{Lemma 2}For every $p>1$ there exist two constants $\kappa_{p}$ and
$K_{p}$, with $0<\kappa_{p}\leq1\leq K_{p}$, such that the following
inequalities hold:%
\begin{gather}
\int_{0}^{1}\left(  \mu^{2}+\left\vert x+ty\right\vert ^{2}\right)
^{\frac{p-2}{2}}(1-t)\,dt\geq\kappa_{p}\left(  \mu^{2}+\left\vert x\right\vert
^{2}+\left\vert y\right\vert ^{2}\right)  ^{\frac{p-2}{2}},\label{prima}\\
\int_{0}^{1}\left(  \mu^{2}+\left\vert x+ty\right\vert ^{2}\right)
^{\frac{p-2}{2}}dt\leq K_{p}\left(  \mu^{2}+\left\vert x\right\vert
^{2}+\left\vert y\right\vert ^{2}\right)  ^{\frac{p-2}{2}} \label{seconda}%
\end{gather}
for every $x,$ $y\in X$ and every constant $\mu\geq0.$
\end{lemma}

\begin{proof}
Let us prove (\ref{prima}). If $1<p\leq2$, then (\ref{prima}) holds with
$\kappa_{p}=2^{\frac{p}{2}-2}$, since
\[
\left(  \mu^{2}+\left\vert x+ty\right\vert ^{2}\right)  ^{\frac{p-2}{2}}%
\geq2^{\frac{p-2}{2}}\left(  \mu^{2}+\left\vert x\right\vert ^{2}+\left\vert
y\right\vert ^{2}\right)  ^{\frac{p-2}{2}}.
\]

If $p>2$, then we consider first the case where $\left\vert y\right\vert
^{2}\leq4\left(  \mu^{2}+\left\vert x\right\vert ^{2}\right)  $. For $0\leq
t\leq1/2$ we have%
\begin{align*}
\left(  \mu^{2}+\left\vert x+ty\right\vert ^{2}\right)  ^{\frac{1}{2}}  &
\geq\left(  \mu^{2}+\left\vert x\right\vert ^{2}\right)  ^{\frac{1}{2}%
}-t\left\vert y\right\vert \geq\left(  1-2t\right)  \left(  \mu^{2}+\left\vert
x\right\vert ^{2}\right)  ^{\frac{1}{2}}\\
&  \geq5^{-\frac{1}{2}}\left(  1-2t\right)  \left(  \mu^{2}+\left\vert
x\right\vert ^{2}+\left\vert y\right\vert ^{2}\right)  ^{\frac{1}{2}},
\end{align*}
hence
\[
\left(  \mu^{2}+\left\vert x+ty\right\vert ^{2}\right)  ^{\frac{p-2}{2}}%
\geq5^{\frac{2-p}{2}}\left(  1-2t\right)  ^{p-2}\left(  \mu^{2}+\left\vert
x\right\vert ^{2}+\left\vert y\right\vert ^{2}\right)  ^{\frac{p-2}{2}}.
\]
We deduce that (\ref{prima}) holds provided
\[
\kappa_{p}\leq5^{\frac{2-p}{2}}\int_{0}^{\frac{1}{2}}\left(  1-2t\right)
^{p-2}(1-t)\,dt=\frac{5^{\frac{2-p}{2}}}{4}\frac{2p-1}{p\left(  p-1\right)
}.
\]

If $p>2$ and $\left\vert y\right\vert ^{2}>4\left(  \mu^{2}+\left\vert
x\right\vert ^{2}\right)  $, then for $1/2\leq t\leq1$ we have%
\begin{align*}
\left(  \mu^{2}+\left\vert x+ty\right\vert ^{2}\right)  ^{\frac{1}{2}}  &
\geq t\left\vert y\right\vert -\left(  \mu^{2}+\left\vert x\right\vert
^{2}\right)  ^{\frac{1}{2}}\geq\left(  t-\frac{1}{2}\right)  \left\vert
y\right\vert \\
&  \geq2\cdot5^{-\frac{1}{2}}\left(  t-\frac{1}{2}\right)  \left(  \mu
^{2}+\left\vert x\right\vert ^{2}+\left\vert y\right\vert ^{2}\right)
^{\frac{1}{2}}.
\end{align*}
Therefore%
\[
\left(  \mu^{2}+\left\vert x+ty\right\vert ^{2}\right)  ^{\frac{p-2}{2}}%
\geq5^{\frac{2-p}{2}}\left(  2t-1\right)  ^{p-2}\left(  \mu^{2}+\left\vert
x\right\vert ^{2}+\left\vert y\right\vert ^{2}\right)  ^{\frac{p-2}{2}}.
\]
We conclude that (\ref{prima}) is verified with
\[
\kappa_{p}\leq5^{\frac{2-p}{2}}\int_{\frac{1}{2}}^{1}\left(  2t-1\right)
^{p-2}(1-t)\,dt=\frac{5^{\frac{2-p}{2}}}{4}\frac{1}{p\left(  p-1\right)  }.
\]
This concludes the proof of (\ref{prima}).

Let us prove (\ref{seconda}). If $p\geq2$ then (\ref{seconda}) holds with
$K_{p}=2^{\frac{p-2}{2}}$, since
\[
\left(  \mu^{2}+\left\vert x+ty\right\vert ^{2}\right)  ^{\frac{p-2}{2}}%
\leq2^{\frac{p-2}{2}}\left(  \mu^{2}+\left\vert x\right\vert ^{2}+\left\vert
y\right\vert ^{2}\right)  ^{\frac{p-2}{2}}.
\]

In the case $1<p<2$ we observe that
\[
\left(  \mu^{2}+\left\vert x+ty\right\vert ^{2}\right)  ^{\frac{p-2}{2}}%
\leq\left\vert \left(  \mu^{2}+\left\vert x\right\vert ^{2}\right)  ^{\frac
{1}{2}}-t\left\vert y\right\vert \right\vert ^{p-2}.
\]
Let $a:=\left(  \mu^{2}+\left\vert x\right\vert ^{2}\right)  ^{\frac{1}{2}}$
and $b:=\left\vert y\right\vert $. Then%
\[
\int_{0}^{1}\left(  \mu^{2}+\left\vert x+ty\right\vert ^{2}\right)
^{\frac{p-2}{2}}dt\leq\int_{0}^{1}\left\vert a-tb\right\vert ^{p-2}dt
\]
and $\left(  \mu^{2}+\left\vert x\right\vert ^{2}+\left\vert y\right\vert
^{2}\right)  ^{\frac{p-2}{2}}=\left(  a^{2}+b^{2}\right)  ^{\frac{p-2}{2}}$.

If $b\leq a$, then%
\[
\int_{0}^{1}\left\vert a-tb\right\vert ^{p-2}dt=\frac{a^{p-1}-\left(
a-b\right)  ^{p-1}}{\left(  p-1\right)  b}\leq\frac{a^{p-2}}{p-1},
\]
where the last inequality is obtained by comparing the difference quotients of
the concave function $t\mapsto t^{p-1}$ in the intervals $[a-b,a]$ and
$[0,a]$. Therefore we have that
\[
\int_{0}^{1}\left\vert a-tb\right\vert ^{p-2}dt\leq\frac{a^{p-2}}{p-1}%
\leq\frac{2^{\frac{2-p}{2}}}{p-1}\left(  a^{2}+b^{2}\right)  ^{\frac{p-2}{2}%
},
\]
and (\ref{seconda}) is satisfied for $K_{p}\geq2^{\frac{2-p}{2}}/(p-1)$.

If $a<b$, then%
\[
\int_{0}^{1}\left\vert a-tb\right\vert ^{p-2}dt=\frac{a^{p-1}+\left(
b-a\right)  ^{p-1}}{\left(  p-1\right)  b}\leq\frac{2^{2-p}b^{p-2}}{p-1}.
\]
On the other hand, in this case we have $b^{p-2}\leq2^{\frac{2-p}{2}}\left(
a^{2}+b^{2}\right)  ^{\frac{p-2}{2}}$, so that (\ref{seconda}) holds for
$K_{p}\geq2^{\frac{3}{2}\left(  2-p\right)  }/(p-1).$
\end{proof}

\begin{lemma}
\label{lemma 8}For every $p>1$, there exist two constants $\theta_{p\text{ }%
}>0$ and $\Theta_{p}>0$ such that for every $\mu\geq0$\ the function $g$
defined in $\left(  \ref{450}\right)  $ satisfies the following inequalities%
\begin{align*}
\theta_{p}\left(  \mu^{2}+\left\vert x\right\vert ^{2}+ \left\vert
y\right\vert ^{2}\right)  ^{\frac{p-2}{2}}\left\vert y\right\vert ^{2}  &
\leq g\left(  x+y\right)  - g\left(  x\right)  -\nabla g\left(  x\right)
\cdot y\\
&  \leq\Theta_{p}\left(  \mu^{2}+\left\vert x\right\vert ^{2}+\left\vert
y\right\vert ^{2}\right)  ^{\frac{p-2}{2}}\left\vert y\right\vert ^{2}%
\end{align*}
for every $x$, $y\in X$.
\end{lemma}

\begin{proof}
By continuity it is enough to prove the statement when $0$ does not belong to
the segment joining $x$ and $x+y.$ In this case the function
\[
h\left(  t\right)  :=\left(  \mu^{2}+\left\vert x+ty\right\vert ^{2}\right)
^{\frac{p}{2}}%
\]
belongs to $C^{\infty}\left(  \left[  0,1\right]  \right)  $, and Taylor's
formula with integral remainder yields
\[
h\left(  1\right)  -h\left(  0\right)  -h^{\prime}\left(  0\right)  =\int
_{0}^{1}h^{\prime\prime}\left(  t\right)  \left(  1-t\right)  \,dt.
\]
By direct computation we see that%
\begin{align*}
p\left(  \left(  p-1\right)  \wedge1\right)   &  \left(  \mu^{2}+\left\vert
x+ty\right\vert ^{2}\right)  ^{\frac{p-2}{2}}\left\vert y\right\vert ^{2} \leq
h^{\prime\prime}\left(  t\right) \\
&  \leq p\left(  \left(  p-1\right)  \vee1\right)  \left(  \mu^{2}+\left\vert
x+ty\right\vert ^{2}\right)  ^{\frac{p-2}{2}}\left\vert y\right\vert ^{2}.
\end{align*}
The conclusion follows from Lemma \ref{Lemma 2}.
\end{proof}

In the proof of Theorem \ref{theorem 2} we will need the following extension
of Lemma \ref{lemma 8} to the family of functions
\[
g_{\beta}\left(  x,y\right)  :=\left(  \mu^{2}+\left\vert x\right\vert
^{2}+\beta^{2}\left\vert y\right\vert ^{2}\right)  ^{\frac{p}{2}},\qquad
\beta\geq0,
\]
defined on the product of two Hilbert spaces $X$ and $Y$.

\begin{lemma}
\label{lemma 9}Let $p>1$, $\beta\geq0$ and $\mu\geq0$. Then%
\begin{align}
g_{\beta}  &  \left(  x+\xi,y+\eta\right)  -g_{\beta}\left(  x,y\right)
-\nabla_{\!x}g_{\beta}\left(  x,y\right)  \cdot\xi-\nabla_{\!y}g_{\beta
}\left(  x,y\right)  \cdot\eta\nonumber\\
\geq &  \theta_{p}\left(  \mu^{2}+\left\vert x\right\vert ^{2}+\left\vert
\xi\right\vert ^{2}+\beta^{2}\left\vert y\right\vert ^{2}+\beta^{2}\left\vert
\eta\right\vert ^{2}\right)  ^{\frac{p-2}{2}}\left(  \left\vert \xi\right\vert
^{2}+\beta^{2}\left\vert \eta\right\vert ^{2}\right) \nonumber
\end{align}
for every $x$, $\xi\in X$, $y$, $\eta\in Y$, where $\theta_{p}$ is the first
constant in Lemma \ref{lemma 8}. Therefore, if $p\geq2$, we have
\begin{align*}
&  g_{\beta}\left(  x+\xi,y+\eta\right)  -g_{\beta}\left(  x,y\right)
-\nabla_{\!x}g_{\beta}\left(  x,y\right)  \cdot\xi-\nabla_{\!y}g_{\beta
}\left(  x,y\right)  \cdot\eta\\
&  \geq\theta_{p}\left(  \mu^{2}+\left\vert x\right\vert ^{2}+\left\vert
\xi\right\vert ^{2}\right)  ^{\frac{p-2}{2}}\left\vert \xi\right\vert
^{2}+\frac{\theta_{p}\beta^{2}}{2}\left(  \mu^{2}+\left\vert x\right\vert
^{2}\right)  ^{\frac{p-2}{2}}\left\vert \eta\right\vert ^{2}+\frac{\theta
_{p}\beta^{p}}{2}\left\vert \eta\right\vert ^{p}%
\end{align*}
for every $x$, $\xi\in X$, $y$, $\eta\in Y$.
\end{lemma}

\begin{proof}
Observing that $g_{\beta}\left(  x,y\right)  =g_{1}\left(  x,\beta y\right)
$, the inequality can be obtained by applying Lemma \ref{lemma 8} to the
Hilbert space $X\times Y$.
\end{proof}

We continue with some technical lemmas which are used in the proofs of
Theorems \ref{theorem 2} and \ref{theorem 1}.

\begin{lemma}
\label{lemma10}Let $X$ be a Hilbert space and let $1<p\leq2$. Then for every
$\mu\geq0$ and every $0<\varepsilon<1$ we have%
\begin{align*}
\left(  \mu^{2}+\left\vert x+y\right\vert ^{2}\right)  ^{\frac{p-2}{2}%
}\left\vert x+y\right\vert ^{2}  &  \leq2\left(  \mu^{2}+\left\vert
x\right\vert ^{2}\right)  ^{\frac{p-2}{2}}\left\vert x\right\vert
^{2}+2\left(  \mu^{2}+\left\vert y\right\vert ^{2}\right)  ^{\frac{p-2}{2}%
}\left\vert y\right\vert ^{2},\\
\varepsilon^{\frac{2-p}{2}}\left(  \mu^{2}+\left\vert y\right\vert
^{2}\right)  ^{\frac{p-2}{2}}\left\vert y\right\vert ^{2}  &  \leq\left(
\mu^{2}+\left\vert x\right\vert ^{2}+\left\vert y\right\vert ^{2}\right)
^{\frac{p-2}{2}}\left\vert y\right\vert ^{2}+\varepsilon\left(  \mu
^{2}+\left\vert x\right\vert ^{2}\right)  ^{\frac{p-2}{2}}\left\vert
x\right\vert ^{2}%
\end{align*}
for every $x,$ $y\in X$.
\end{lemma}

\begin{proof}
Since the mapping $t\mapsto\left(  \mu^{2}+t\right)  ^{\frac{p-2}{2}}t$ is
nondecreasing, while the mapping $t\mapsto\left(  \mu^{2}+t\right)
^{\frac{p-2}{2}}$ is nonincreasing, we have%
\begin{align*}
\left(  \mu^{2}+\left\vert x+y\right\vert ^{2}\right)  ^{\frac{p-2}{2}}  &
\left\vert x+y\right\vert ^{2}\leq\left(  \mu^{2}+2\left\vert x\right\vert
^{2}+2\left\vert y\right\vert ^{2}\right)  ^{\frac{p-2}{2}}\left(  2\left\vert
x\right\vert ^{2}+2\left\vert y\right\vert ^{2}\right) \\
&  \leq2\left(  \mu^{2}+\left\vert x\right\vert ^{2}\right)  ^{\frac{p-2}{2}%
}\left\vert x\right\vert ^{2}+2\left(  \mu^{2}+\left\vert y\right\vert
^{2}\right)  ^{\frac{p-2}{2}}\left\vert y\right\vert ^{2},
\end{align*}
which proves the first inequality. For the same reason, we have%
\begin{align*}
\varepsilon^{\frac{2-p}{2}}\left(  \mu^{2}+\left\vert y\right\vert
^{2}\right)  ^{\frac{p-2}{2}}  &  \left\vert y\right\vert ^{2}\leq
\varepsilon^{\frac{2-p}{2}}\left(  \mu^{2}+\varepsilon\left\vert x\right\vert
^{2}+\left\vert y\right\vert ^{2}\right)  ^{\frac{p-2}{2}}\left(
\varepsilon\left\vert x\right\vert ^{2}+\left\vert y\right\vert ^{2}\right) \\
&  \leq\left(  \mu^{2}+\left\vert x\right\vert ^{2}+\left\vert y\right\vert
^{2}\right)  ^{\frac{p-2}{2}}\left(  \varepsilon\left\vert x\right\vert
^{2}+\left\vert y\right\vert ^{2}\right) \\
&  \leq\varepsilon\left(  \mu^{2}+\left\vert x\right\vert ^{2}\right)
^{\frac{p-2}{2}}\left\vert x\right\vert ^{2}+\left(  \mu^{2}+\left\vert
x\right\vert ^{2}+\left\vert y\right\vert ^{2}\right)  ^{\frac{p-2}{2}%
}\left\vert y\right\vert ^{2},
\end{align*}
and this concludes the proof.
\end{proof}

\begin{lemma}
\label{lemma 12} Let $1<p\leq2$. Then
\begin{equation}
b^{p}\leq8 \varepsilon^{\frac{p-2}{p}}\left(  \mu^{2}+a^{2}+b^{2}\right)
^{\frac{p-2}{2}}b^{2}+\varepsilon a^{p}+\varepsilon\mu^{p} \label{400}%
\end{equation}
for every $a\geq0$, $b\geq0$, $\mu\geq0$, and $0<\varepsilon<1$.
\end{lemma}

\begin{proof}
By Young's inequality with exponents $2/p$ and $2/\left(  2-p\right)  $ we
have that%
\begin{align}
b^{p}  &  =\left(  2^{\frac{2-p}{2}} \varepsilon^{\frac{p-2}{2}}\left(
\mu^{2}+b^{2} \right)  ^{\frac{p\left(  p-2\right)  }{4}} b^{p} \right)
\left(  2^{\frac{p-2}{2}} \varepsilon^{\frac{2-p}{2}}\left(  \mu^{2}%
+b^{2}\right)  ^{\frac{p\left(  2-p\right)  }{4}}\right) \label{960}\\
&  \leq2 \varepsilon^{\frac{p-2}{p}}\left(  \mu^{2}+b^{2}\right)  ^{\frac
{p-2}{2}}b^{2}+ \frac{\varepsilon}{2} \left(  \mu^{2}+b^{2}\right)  ^{\frac
{p}{2}}.\nonumber
\end{align}

If $a^{2}\leq\mu^{2}+b^{2}$, we have%
\[
\left(  \mu^{2}+b^{2}\right)  ^{\frac{p-2}{2}}\leq2\left(  \mu^{2}+a^{2}%
+b^{2}\right)  ^{\frac{p-2}{2}},
\]
and so from (\ref{960}) we obtain
\begin{equation}
b^{p}\leq4\varepsilon^{\frac{p-2}{p}}\left(  \mu^{2}+a^{2}+b^{2}\right)
^{\frac{p-2}{2}}b^{2}+\frac{\varepsilon}{2}\left(  \mu^{p}+b^{p}\right)
.\nonumber
\end{equation}
Subtracting $\left(  \varepsilon/2\right)  b^{p}$ to both sides we obtain
(\ref{400}).

On the other hand, if $a^{2}>\mu^{2}+b^{2}$, then from (\ref{960}) with
$\left(  \mu^{2}+b^{2}\right)  $ replaced by $\left(  \mu^{2}+a^{2}%
+b^{2}\right)  $ we obtain%
\begin{align*}
b^{p}  &  \leq2\varepsilon^{\frac{p-2}{p}}\left(  \mu^{2}+a^{2}+b^{2}\right)
^{\frac{p-2}{2}}b^{2}+\frac{\varepsilon}{2}\left(  \mu^{2}+a^{2}+b^{2}\right)
^{\frac{p}{2}}\\
&  \leq2\varepsilon^{\frac{p-2}{p}}\left(  \mu^{2}+a^{2}+b^{2}\right)
^{\frac{p-2}{2}}b^{2}+\varepsilon a^{p},
\end{align*}
which proves (\ref{400}).
\end{proof}

The next lemma shows that condition (\ref{3}) in Theorem \ref{theorem 2} can
be obtained from a suitable bound on the second derivatives of $f$. This is
trivial in the case $p\geq2$, but requires some work in the case $1<p<2$.

\begin{lemma}
\label{Lemma 3}Let $X$ be a Hilbert space, and let $f\in C^{1}\left(
X\right)  \cap C^{2}\left(  X\setminus\left\{  0\right\}  \right)  .$ Assume
that there exist two constants $p>1$, $C>0$, and $\mu\geq0$ such that%
\begin{equation}
\left\vert \nabla^{2}f\left(  x\right)  \right\vert \leq C\left(  \mu
^{2}+\left\vert x\right\vert ^{2}\right)  ^{\frac{p-2}{2}} \label{12}%
\end{equation}
for every $x\in X\setminus\left\{  0\right\}  $. Then
\begin{equation}
\left\vert \nabla f\left(  x+y\right)  -\nabla f\left(  x\right)  \right\vert
\leq K_{p}C\left(  \mu^{2}+\left\vert x\right\vert ^{2}+\left\vert
y\right\vert ^{2}\right)  ^{\frac{p-2}{2}}\left\vert y\right\vert \label{13}%
\end{equation}
for every $x$, $y\in X$, where $K_{p}$ is the second constant in Lemma
\ref{Lemma 2}.
\end{lemma}

\begin{proof}
By continuity it is enough to prove the statement when $0$ does not belong to
the segment joining $x$ and $x+y.$ In this case by (\ref{12}) we have%
\[
\left\vert \nabla f\left(  x+y\right)  -\nabla f\left(  x\right)  \right\vert
\leq C\left\vert y\right\vert \int_{0}^{1}\left(  \mu^{2}+\left\vert
x+ty\right\vert ^{2}\right)  ^{\frac{p-2}{2}}dt,
\]
and the conclusion follows from Lemma\nolinebreak\ \ref{Lemma 2}.
\end{proof}

The estimate given by the following lemma will be crucial in the proof of
Theorem \ref{theorem 2}.

\begin{lemma}
\label{Lemma 4}Let $X$ be a Hilbert space and let $f\in C^{1}\left(  X\right)
$. Assume that there exist $p>1$ and $\mu\geq0$ such that
\begin{equation}
\left\vert \nabla f\left(  x+y\right)  -\nabla f\left(  x\right)  \right\vert
\leq\left(  \mu^{2}+\left\vert x\right\vert ^{2}+\left\vert y\right\vert
^{2}\right)  ^{\frac{p-2}{2}}\left\vert y\right\vert \label{50}%
\end{equation}
for every $x$, $y\in X$. If $1<p\leq2$, then for every $\varepsilon>0$ there
exists $c_{\varepsilon,p}>0$, depending only on $\varepsilon$ and $p$, such
that
\begin{align}
\left\vert f\left(  x+y+z\right)  -f\left(  x+y\right)  -\nabla f\left(
x\right)  \cdot z\right\vert  & \label{14}\\
\leq\varepsilon\left(  \mu^{2}+\left\vert x\right\vert ^{2}+\left\vert
y\right\vert ^{2}\right)  ^{\frac{p-2}{2}}  &  \left\vert y\right\vert
^{2}+c_{\varepsilon,p}\left(  \mu^{2}+\left\vert z\right\vert ^{2}\right)
^{\frac{p-2}{2}}\left\vert z\right\vert ^{2}\nonumber
\end{align}
for every $x$, $y$, $z\in X$.

If $p\geq2$, then for every $\varepsilon>0$ there exists $c_{\varepsilon,p}%
>0$, depending only on $\varepsilon$ and $p$, such that%
\begin{align}
\left\vert f\left(  x+y+z\right)  -f\left(  x+y\right)  -\nabla f\left(
x\right)  \cdot z\right\vert  & \label{61}\\
\leq\varepsilon\left(  \mu^{2}+\left\vert x\right\vert ^{2}+\left\vert
y\right\vert ^{2}\right)  ^{\frac{p-2}{2}}  &  \left\vert y\right\vert
^{2}+c_{\varepsilon,p}\left(  \mu^{2}+\left\vert x\right\vert ^{2}\right)
^{\frac{p-2}{2}}\left\vert z\right\vert ^{2}+c_{\varepsilon,p}\left\vert
z\right\vert ^{p}\nonumber
\end{align}
for every $x$, $y$, $z\in X$.
\end{lemma}

\begin{proof}
Let us consider first the case $1<p\leq2$. Clearly (\ref{50}) implies that
\[
\left|  \nabla f\left(  x+y\right)  -\nabla f\left(  x\right)  \right|
\leq\left(  \mu^{2}+\left|  y\right|  ^{2}\right)  ^{\frac{p-2}{2}}\left|
y\right|
\]
for every $x$, $y\in X$. By the Mean Value Theorem we have
\begin{align*}
\left|  f\left(  x+y\right)  -f\left(  x\right)  -\nabla f\left(  x\right)
\cdot y\right|   &  \leq\left|  \nabla f\left(  x+ty\right)  -\nabla f\left(
x\right)  \right|  \left|  y\right| \\
&  \leq\left(  \mu^{2}+t^{2}\left|  y\right|  ^{2}\right)  ^{\frac{p-2}{2}%
}t\left|  y\right|  ^{2}%
\end{align*}
for some $t\in\left[  0,1\right]  $. Since the function $t\mapsto\left(
\mu^{2}+t^{2}\left|  y\right|  ^{2}\right)  ^{\frac{p-2}{2}}t$ is
nondecreasing, we obtain
\begin{equation}
\left|  f\left(  x+y\right)  -f\left(  x\right)  -\nabla f\left(  x\right)
\cdot y\right|  \leq\left(  \mu^{2}+\left|  y\right|  ^{2}\right)
^{\frac{p-2}{2}}\left|  y\right|  ^{2} \label{51}%
\end{equation}
for every $x$, $y\in X$.

By (\ref{50}) and (\ref{51}) we have%
\begin{align}
|f(x+y  &  +z)-f\left(  x+y\right)  -\nabla f\left(  x\right)  \cdot
z|\nonumber\\
\leq &  \left|  f\left(  x+y+z\right)  -f\left(  x+y\right)  -\nabla f\left(
x+y\right)  \cdot z\right| \label{52}\\
&  +\left|  \nabla f\left(  x+y\right)  \cdot z-\nabla f\left(  x\right)
\cdot z\right| \nonumber\\
\leq &  \left(  \mu^{2}+\left|  z\right|  ^{2}\right)  ^{\frac{p-2}{2}}\left|
z\right|  ^{2}+\left(  \mu^{2}+\left|  x\right|  ^{2}+\left|  y\right|
^{2}\right)  ^{\frac{p-2}{2}}\left|  y\right|  \left|  z\right|  .\nonumber
\end{align}
We now estimate the last term. If $\left|  z\right|  \leq\mu$ then for every
$\varepsilon>0$ we have%
\begin{align}
\left(  \mu^{2}+\left|  x\right|  ^{2}+\left|  y\right|  ^{2}\right)
^{\frac{p-2}{2}}  &  \left|  y\right|  \left|  z\right|  \leq\varepsilon
\left(  \mu^{2}+\left|  x\right|  ^{2}+\left|  y\right|  ^{2}\right)
^{\frac{p-2}{2}}\left|  y\right|  ^{2}\nonumber\\
&  +\frac{1}{4\varepsilon}\left(  \mu^{2}+\left|  x\right|  ^{2}+\left|
y\right|  ^{2}\right)  ^{\frac{p-2}{2}}\left|  z\right|  ^{2}\label{53}\\
\leq &  \varepsilon\left(  \mu^{2}+\left|  x\right|  ^{2}+\left|  y\right|
^{2}\right)  ^{\frac{p-2}{2}}\left|  y\right|  ^{2}+\frac{1}{4\varepsilon}%
\mu^{p-2}\left|  z\right|  ^{2}\nonumber\\
\leq &  \varepsilon\left(  \mu^{2}+\left|  x\right|  ^{2}+\left|  y\right|
^{2}\right)  ^{\frac{p-2}{2}}\left|  y\right|  ^{2}+\frac{1}{\varepsilon
}\left(  \mu^{2}+\left|  z\right|  ^{2}\right)  ^{\frac{p-2}{2}}\left|
z\right|  ^{2}.\nonumber
\end{align}
If $\left|  z\right|  >\mu$, let $q:=p/\left(  p-1\right)  $ be the conjugate
exponent of $p$. Since $\frac{p-2}{2}+\frac{1}{2}-\frac{1}{q}=\frac{p-2}{2q}$
we have
\begin{align}
\left(  \mu^{2}+\left|  x\right|  ^{2}+\left|  y\right|  ^{2}\right)
^{\frac{p-2}{2}}  &  \left|  y\right|  \left|  z\right|  =\left(  \mu
^{2}+\left|  x\right|  ^{2}+\left|  y\right|  ^{2}\right)  ^{\frac{p-2}{2}%
}\left|  y\right|  ^{1-\frac{2}{q}}\left|  y\right|  ^{\frac{2}{q}}\left|
z\right| \nonumber\\
&  \leq\left(  \mu^{2}+\left|  x\right|  ^{2}+\left|  y\right|  ^{2}\right)
^{\frac{p-2}{2q}}\left|  y\right|  ^{\frac{2}{q}}\left|  z\right| \label{54}\\
&  \leq\varepsilon\left(  \mu^{2}+\left|  x\right|  ^{2}+\left|  y\right|
^{2}\right)  ^{\frac{p-2}{2}}\left|  y\right|  ^{2}+\frac{1}{p\left(
q\varepsilon\right)  ^{p-1}}\left|  z\right|  ^{p}\nonumber\\
&  \leq\varepsilon\left(  \mu^{2}+\left|  x\right|  ^{2}+\left|  y\right|
^{2}\right)  ^{\frac{p-2}{2}}\left|  y\right|  ^{2}+k_{\varepsilon,p}\left(
\mu^{2}+\left|  z\right|  ^{2}\right)  ^{\frac{p-2}{2}}\left|  z\right|
^{2}\nonumber
\end{align}
for some constant $k_{\varepsilon,p}$ depending only on $\varepsilon$ and $p$.

Let us consider now the case $p\geq2$. By the Mean Value Theorem and by the
Cauchy Inequality we have%
\begin{align*}
|f(x  &  +y+z)-f\left(  x+y\right)  -\nabla f\left(  x\right)  \cdot z|\\
\leq &  \left\vert f\left(  x+y+z\right)  -f\left(  x+y\right)  -\nabla
f\left(  x+y\right)  \cdot z\right\vert \\
&  +\left\vert \nabla f\left(  x+y\right)  \cdot z-\nabla f\left(  x\right)
\cdot z\right\vert \\
\leq &  \left(  \mu^{2}+2\left\vert x\right\vert ^{2}+2\left\vert y\right\vert
^{2}+\left\vert z\right\vert ^{2}\right)  ^{\frac{p-2}{2}}\left\vert
z\right\vert ^{2}+\left(  \mu^{2}+\left\vert x\right\vert ^{2}+\left\vert
y\right\vert ^{2}\right)  ^{\frac{p-2}{2}}\left\vert y\right\vert \left\vert
z\right\vert \\
\leq &  6^{\frac{p-2}{2}}\left(  \mu^{2}+\left\vert x\right\vert ^{2}\right)
^{\frac{p-2}{2}}\left\vert z\right\vert ^{2}+6^{\frac{p-2}{2}}\left\vert
y\right\vert ^{p-2}\left\vert z\right\vert ^{2}+3^{\frac{p-2}{2}}\left\vert
z\right\vert ^{p}\\
&  +\frac{\varepsilon}{2}\left(  \mu^{2}+\left\vert x\right\vert
^{2}+\left\vert y\right\vert ^{2}\right)  ^{\frac{p-2}{2}}\left\vert
y\right\vert ^{2}+k_{\varepsilon,p}\left(  \mu^{2}+\left\vert x\right\vert
^{2}\right)  ^{\frac{p-2}{2}}\left\vert z\right\vert ^{2}+k_{\varepsilon
,p}\left\vert y\right\vert ^{p-2}\left\vert z\right\vert ^{2}%
\end{align*}
for some constant $k_{\varepsilon,p}$ depending only on $\varepsilon$ and $p$.
The conclusion follows from Young' inequality with exponents $p/\left(
p-2\right)  $ and $p/2$.
\end{proof}

The following lemma states an elementary property of strictly $2$-quasiconvex functions.

\begin{lemma}
\label{lemma 11}Assume that $f:\mathbb{M}_{\operatorname*{sym}}^{n\times
n}\rightarrow\mathbb{R}$ satisfies the strict $2$-quasiconvexity condition (a)
of Theorem \ref{theorem 2} for some constants $p>1$, $\mu\geq0$, $\nu>0$, and
let $g:\mathbb{M}_{\operatorname*{sym}}^{n\times n}\rightarrow\mathbb{R}$ be
the function defined by%
\begin{equation}
g(A):=\left(  \mu^{2}+\left\vert A\right\vert ^{2}\right)  ^{\frac{p}{2}}.
\label{950}%
\end{equation}
Then the function $f_{\lambda}:=f-\lambda g$ is $2$-quasiconvex for
$\lambda\leq\nu/\Theta_{p}$, where $\Theta_{p}$ is the second constant in
Lemma \ref{lemma 8}.
\end{lemma}

\begin{proof}
Let $A\in\mathbb{M}_{\operatorname*{sym}}^{n\times n}$ and $\phi\in C_{c}%
^{2}\left(  Q\right)  $. Since, by periodicity,%
\[
\int_{Q}\nabla g(A)\cdot\nabla^{2}\phi\,dx=0,
\]
we have%
\begin{align*}
\int_{Q}  &  \left[  f_{\lambda}\left(  A+\nabla^{2}\phi\right)  -f_{\lambda
}\left(  A\right)  \right]  \,dx=\int_{Q}\left[  f\left(  A+\nabla^{2}%
\phi\right)  -f\left(  A\right)  \right]  \,dx\\
&  -\lambda\int_{Q}\left[  g\left(  A+\nabla^{2}\phi\right)  -g\left(
A\right)  +\nabla g(A)\cdot\nabla^{2}\phi\right]  \,dx\\
&  \geq\left(  \nu-\lambda\Theta_{p}\right)  \int_{Q}\left(  \mu^{2}+\left|
A\right|  ^{2}+\left|  \nabla^{2}\phi\right|  ^{2}\right)  ^{\frac{p-2}{2}%
}\left|  \nabla^{2}\phi\right|  ^{2}\,dx\geq0,
\end{align*}
which concludes the proof.
\end{proof}

In the proof of Theorem \ref{theorem 2} we need the following generalization
of Lemma \ref{lemma 1}.

\begin{lemma}
\label{lemma 5}For every $p>1$ there exists a constant $\tau_{n,p}\geq1$ such
that
\begin{equation}
\int_{Q}\left(  \mu^{2}+\left\vert \nabla\psi^{s}\right\vert ^{2}\right)
^{\frac{p-2}{2}}\left\vert \nabla\psi^{s}\right\vert ^{2}\,dx\leq\tau
_{n,p}\int_{Q}\left(  \mu^{2}+\left\vert \nabla\psi^{a}\right\vert
^{2}\right)  ^{\frac{p-2}{2}}\left\vert \nabla\psi^{a}\right\vert ^{2}\,dx
\label{16}%
\end{equation}
for every constant $\mu\geq0$ and every $Q$-periodic function $\psi
:\mathbb{R}^{n}\rightarrow\mathbb{R}^{n}$ of class $C^{\infty}$ with
$\operatorname*{div}\psi=0$.
\end{lemma}

\begin{proof}
Let $\mu$ and $\psi$ be as in the statement of the lemma. In the case $p\geq2$
by Lemma \ref{lemma 1} we have%
\[
\int_{Q}\left\vert \nabla\psi\right\vert ^{p}\,dx\leq\gamma_{n,p}\int
_{Q}\left\vert \nabla\psi^{a}\right\vert ^{p}\,dx,
\]
and so
\begin{align*}
\int_{Q}  &  \left(  \mu^{2}+\left\vert \nabla\psi\right\vert ^{2}\right)
^{\frac{p-2}{2}}\left\vert \nabla\psi\right\vert ^{2}\,dx\\
&  \leq2^{\frac{p-2}{2}}\left\{  \mu^{p-2}\int_{Q}\left\vert \nabla
\psi\right\vert ^{2}\,dx+\int_{Q}\left\vert \nabla\psi\right\vert
^{p}\,dx\right\} \\
&  \leq2^{\frac{p-2}{2}}\left\{  \mu^{p-2}\gamma_{n,2}\int_{Q}\left\vert
\nabla\psi^{a}\right\vert ^{2}\,dx+\gamma_{n,p}\int_{Q}\left\vert \nabla
\psi^{a}\right\vert ^{p}\,dx\right\} \\
&  \leq\tau_{n,p}\int_{Q}\left(  \mu^{2}+\left\vert \nabla\psi^{a}\right\vert
^{2}\right)  ^{\frac{p-2}{2}}\left\vert \nabla\psi^{a}\right\vert ^{2}\,dx
\end{align*}
with $\tau_{n,p}:=2^{\frac{p-2}{2}}\left(  \gamma_{n,2}+\gamma_{n,p},\right)
$.

We consider now the case $1<p<2$. Let $E_{\mu}:=\left\{  \left\vert \nabla
\psi^{a}\right\vert \leq\mu\right\}  $ and $E^{\mu}:=\left\{  \left\vert
\nabla\psi^{a}\right\vert >\mu\right\}  $, and let $\Psi_{\mu}:=1_{E_{\mu}%
}\nabla\psi^{a}$ and $\Psi^{\mu}:=1_{E^{\mu}}\nabla\psi^{a}$, where $1_{E}$ is
the characteristic function of the set $E$. Note that $\Psi_{\mu}$ and
$\Psi^{\mu}$ are periodic vector-fields. Let $\psi_{\mu}$ and $\psi^{\mu}$ be
periodic solutions of the equations
\[
\Delta\psi_{\mu}=2\operatorname*{div}\Psi_{\mu}\qquad\text{and}\qquad
\Delta\psi^{\mu}=2\operatorname*{div}\Psi^{\mu}.
\]
{}From the first equation we get
\begin{equation}
\int_{Q}\left\vert \nabla\psi_{\mu}\right\vert ^{2}\,dx\leq4\int_{Q}\left\vert
\Psi_{\mu}\right\vert ^{2}\,dx. \label{30}%
\end{equation}
Standard $L^{p}$ estimates for periodic solutions of the Poisson equation (see
\cite{GT}) yield a constant $\tilde{\gamma}_{n,p}\geq4$ such that
\begin{equation}
\int_{Q}\left\vert \nabla\psi^{\mu}\right\vert ^{p}\,dx\leq\tilde{\gamma
}_{n,p}\int_{Q}\left\vert \Psi^{\mu}\right\vert ^{p}\,dx. \label{31}%
\end{equation}
{}From (\ref{30}) we obtain
\begin{align}
\int_{Q}\left(  \mu^{2}+\left\vert \nabla\psi_{\mu}\right\vert ^{2}\right)
^{\frac{p-2}{2}}  &  \left\vert \nabla\psi_{\mu}\right\vert ^{2}\,dx\leq
\mu^{p-2}\int_{Q}\left\vert \nabla\psi_{\mu}\right\vert ^{2}\,dx\nonumber\\
&  \leq4\mu^{p-2}\int_{Q}\left\vert \Psi_{\mu}\right\vert ^{2}\,dx\label{40}\\
&  \leq8\int_{Q\cap E_{\mu}}\left(  \mu^{2}+\left\vert \nabla\psi
^{a}\right\vert ^{2}\right)  ^{\frac{p-2}{2}}\left\vert \nabla\psi
^{a}\right\vert ^{2}\,dx,\nonumber
\end{align}
where the last inequality follows from the fact that $\left\vert \Psi_{\mu
}\right\vert =\left\vert \nabla\psi^{a}\right\vert 1_{E_{\mu}}\leq\mu$. {}From
(\ref{31}) we obtain
\begin{align}
\int_{Q}\left(  \mu^{2}+\left\vert \nabla\psi^{\mu}\right\vert ^{2}\right)
^{\frac{p-2}{2}}  &  \left\vert \nabla\psi^{\mu}\right\vert ^{2}\,dx\leq
\int_{Q}\left\vert \nabla\psi^{\mu}\right\vert ^{p}\,dx\nonumber\\
&  \leq\tilde{\gamma}_{n,p}\int_{Q}\left\vert \Psi^{\mu}\right\vert
^{p}\,dx\label{41}\\
&  \leq2\tilde{\gamma}_{n,p}\int_{Q\cap E^{_{\mu}}}\left(  \mu^{2}+\left\vert
\nabla\psi^{a}\right\vert ^{2}\right)  ^{\frac{p-2}{2}}\left\vert \nabla
\psi^{a}\right\vert ^{2}\,dx,\nonumber
\end{align}
where the last inequality follows from the fact that $\left\vert \Psi^{_{\mu}%
}\right\vert =\left\vert \nabla\psi^{a}\right\vert 1_{E^{_{\mu}}}\geq
\mu1_{E^{_{\mu}}}$. By (\ref{901}) we have
\[
\Delta\left(  \psi_{\mu}+\psi^{\mu}\right)  =2\operatorname*{div}\nabla
\psi^{a}=\Delta\psi,
\]
and since $\psi_{\mu}+\psi^{\mu}-\psi$ is a periodic function we deduce that
$\nabla\psi=\nabla\psi_{\mu}+\nabla\psi^{\mu}$. Finally, from Lemma
\ref{lemma10} and using (\ref{40}) and (\ref{41}) we obtain%
\[
\int_{Q}\left(  \mu^{2}+\left\vert \nabla\psi\right\vert ^{2}\right)
^{\frac{p-2}{2}}\left\vert \nabla\psi\right\vert ^{2}\,dx\leq\tau_{n,p}%
\int_{Q}\left(  \mu^{2}+\left\vert \nabla\psi^{a}\right\vert ^{2}\right)
^{\frac{p-2}{2}}\left\vert \nabla\psi^{a}\right\vert ^{2}\,dx,
\]
with $\tau_{n,p}=2\tilde{\gamma}_{n,p}$. Since $\left\vert \nabla\psi
^{s}\right\vert \leq\left\vert \nabla\psi\right\vert $ and the mapping
$t\mapsto\left(  \mu^{2}+t\right)  ^{\frac{p-2}{2}}t$ is
nondecreasing,\ inequality (\ref{16}) follows.
\end{proof}

\section{Proofs\label{section proofs}}

\begin{proof}
[Proof of Theorem \ref{theorem 2}]We begin by observing that (\ref{3}) gives
\begin{equation}
\left|  f\left(  A\right)  \right|  \leq k_{f}\left(  1+\left|  A\right|
^{p}\right)  \qquad\forall A\in\mathbb{M}_{\operatorname*{sym}}^{n\times n}
\label{18}%
\end{equation}
for a suitable constant $k_{f}$ depending on $f$.

\noindent\textbf{Step 1:} We first consider the case $1<p<2$. Let
$g:\mathbb{M}^{n\times n}\rightarrow\mathbb{R}$ be the function defined by%
\[
g(A):=\left(  \mu^{2}+\left\vert A^{a}\right\vert ^{2}\right)  ^{\frac{p}{2}%
}-\mu^{p}.
\]
Given a constant $\beta>0$, to be chosen at the end of the proof, let
$G:\mathbb{M}^{n\times n}\rightarrow\mathbb{R}$ be the function defined by%
\begin{equation}
G(A):=f(A^{s})+\beta g(A^{a}), \label{19}%
\end{equation}
and let $F$ be its $1$-quasiconvexification, i.e., (see, e.g., \cite{D})
\begin{equation}
F(A)=\inf\left\{  \int_{Q}G(A+\nabla\varphi(x))\,dx:\,\varphi\in
C_{\operatorname*{per}}^{\infty}(Q;\mathbb{R}^{n})\right\}  , \label{15}%
\end{equation}
for all $A\in\mathbb{M}^{n\times n}$.

We want to prove that for every $\varepsilon>0$ there exists $\beta>0$ such
that
\begin{equation}
\int_{Q}\left[  G\left(  A+\nabla\varphi\right)  -G\left(  A\right)  \right]
\,dx\geq-\varepsilon\left(  \mu^{2}+\left|  A^{a}\right|  ^{2}\right)
^{\frac{p-2}{2}}\left|  A^{a}\right|  ^{2} \label{20}%
\end{equation}
for every $A\in\mathbb{M}^{n\times n}$ and for every $\varphi\in
C_{\operatorname*{per}}^{\infty}(Q;\mathbb{R}^{n})$. In view of (\ref{15})
this will imply that for every $A\in\mathbb{M}^{n\times n}$ we have%
\begin{equation}
G(A)-\varepsilon\left(  \mu^{2}+\left|  A^{a}\right|  ^{2}\right)
^{\frac{p-2}{2}}\left|  A^{a}\right|  ^{2}\leq F(A)\leq G(A) \label{23}%
\end{equation}
which yields (\ref{5}). Inequality (\ref{6}) follows from (\ref{18}),
(\ref{19}) and (\ref{23}).

Let us prove (\ref{20}). Fix $\varphi\in C_{\operatorname*{per}}^{\infty
}(Q;\mathbb{R}^{n})$ and consider the periodic Helmholtz decomposition
\[
\varphi=\nabla\phi+\psi
\]
given by Lemma \ref{lemma 6}. Following the argument used by M\"{u}ller and
\v{S}ver\'{a}k in the proof of Lemma 4.2 in \cite{MS}, we have%
\begin{align}
\int_{Q}[G  &  \left(  A+\nabla\varphi\right)  -G\left(  A\right)
]\,dx\nonumber\\
=  &  \int_{Q}\left[  f\left(  A^{s}+\nabla^{2}\phi+\nabla\psi^{s}\right)
-f\left(  A^{s}+\nabla^{2}\phi\right)  \right]  \,dx\nonumber\\
&  +\int_{Q}\left[  f\left(  A^{s}+\nabla^{2}\phi\right)  -f\left(
A^{s}\right)  \right]  \,dx\label{17}\\
&  +\beta\int_{Q}\left[  \left(  \mu^{2}+\left\vert A^{a}+\nabla\psi
^{a}\right\vert ^{2}\right)  ^{\frac{p}{2}}-\left(  \mu^{2}+\left\vert
A^{a}\right\vert ^{2}\right)  ^{\frac{p}{2}}\right]  \,dx\nonumber\\
=  &  \!:I_{1}+I_{2}+I_{3}.\nonumber
\end{align}
Since $\nabla f\left(  A^{s}\right)  $ is a symmetric matrix we have $\nabla
f\left(  A^{s}\right)  \cdot\nabla\psi^{s}=\nabla f\left(  A^{s}\right)
\cdot\nabla\psi$, and therefore, by periodicity,
\[
\int_{Q}\nabla f\left(  A^{s}\right)  \cdot\nabla\psi^{s}\,dx=0.
\]
Hence
\[
I_{1}=\int_{Q}\left[  f\left(  A^{s}+\nabla^{2}\phi+\nabla\psi^{s}\right)
-f\left(  A^{s}+\nabla^{2}\phi\right)  -\nabla f\left(  A^{s}\right)
\cdot\nabla\psi^{s}\right]  \,dx.
\]
By Lemma \ref{Lemma 4} we have%
\begin{align*}
I_{1}  &  \geq-\nu\int_{Q}\left(  \mu^{2}+\left\vert A^{s}\right\vert
^{2}+\left\vert \nabla^{2}\phi\right\vert ^{2}\right)  ^{\frac{p-2}{2}%
}\left\vert \nabla^{2}\phi\right\vert ^{2}\,dx\\
&  -c_{\nu,p,L}\int_{Q}\left(  \mu^{2}+\left\vert \nabla\psi^{s}\right\vert
^{2}\right)  ^{\frac{p-2}{2}}\left\vert \nabla\psi^{s}\right\vert ^{2}\,dx,
\end{align*}
while the strict $2$-quasiconvexity of $f$ (condition (a)) yields%
\[
I_{2}\geq\nu\int_{Q}\left(  \mu^{2}+\left\vert A^{s}\right\vert ^{2}%
+\left\vert \nabla^{2}\phi\right\vert ^{2}\right)  ^{\frac{p-2}{2}}\left\vert
\nabla^{2}\phi\right\vert ^{2}\,dx,
\]
and so, using Lemma \ref{lemma 5},\ we obtain%
\begin{align}
I_{1}+I_{2}  &  \geq-c_{\nu,p,L}\int_{Q}\left(  \mu^{2}+\left\vert \nabla
\psi^{s}\right\vert ^{2}\right)  ^{\frac{p-2}{2}}\left\vert \nabla\psi
^{s}\right\vert ^{2}\,dx\label{24}\\
&  \geq-c_{\nu,p,L}\,\tau_{n,p}\int_{Q}\left(  \mu^{2}+\left\vert \nabla
\psi^{a}\right\vert ^{2}\right)  ^{\frac{p-2}{2}}\left\vert \nabla\psi
^{a}\right\vert ^{2}\,dx.\nonumber
\end{align}
Since $\nabla g(A^{a})$ is an antisymmetric matrix we have $\nabla
g(A^{a})\cdot\nabla\psi^{a}=\nabla g(A^{a})\cdot\nabla\psi$, and therefore, by
periodicity,
\[
\int_{Q}\nabla g(A^{a})\cdot\nabla\psi^{a}\,dx=0.
\]
Hence, by Lemma \ref{lemma 8} and Lemma \ref{lemma10}, for every $0<\delta<1$
we obtain%
\begin{align*}
I_{3}  &  =\beta\int_{Q}\left[  g\left(  A^{a}+\nabla\psi^{a}\right)
-g\left(  A^{a}\right)  -\nabla g(A^{a})\cdot\nabla\psi^{a}\right]  \,dx\\
&  \geq\beta\theta_{p}\int_{Q}\left(  \mu^{2}+\left\vert A^{a}\right\vert
^{2}+\left\vert \nabla\psi^{a}\right\vert ^{2}\right)  ^{\frac{p-2}{2}%
}\left\vert \nabla\psi^{a}\right\vert ^{2}\,dx\\
&  \geq\beta\theta_{p}\delta^{\frac{2-p}{2}}\int_{Q}\left(  \mu^{2}+\left\vert
\nabla\psi^{a}\right\vert ^{2}\right)  ^{\frac{p-2}{2}}\left\vert \nabla
\psi^{a}\right\vert ^{2}\,dx-\beta\theta_{p}\delta\left(  \mu^{2}+\left\vert
A^{a}\right\vert ^{2}\right)  ^{\frac{p-2}{2}}\left\vert A^{a}\right\vert
^{2}.
\end{align*}
Choosing $\beta>0$ and $0<\delta<1$ so that
\[
\beta\theta_{p}\delta^{\frac{2-p}{2}}\geq c_{\nu,p,L}\,\tau_{n,p},\qquad
\beta\theta_{p}\delta\leq\varepsilon
\]
we obtain
\[
I_{1}+I_{2}+I_{3}\geq-\varepsilon\left(  \mu^{2}+\left\vert A^{a}\right\vert
^{2}\right)  ^{\frac{p-2}{2}}\left\vert A^{a}\right\vert ^{2},
\]
which, together with (\ref{17}), yields (\ref{20}).

\noindent\textbf{Step 2:} Let us consider now the case $p\geq2$. Let
$\lambda:=\nu/\Theta_{p}$, where $\Theta_{p}$ is the second constant in Lemma
\ref{lemma 8}. Given a constant $\beta>0$, to be chosen at the end of the
proof, let $F:\mathbb{M}^{n\times n}\rightarrow\mathbb{R}$ be the function
defined by%
\begin{equation}
F(A):=f(A^{s})-\lambda\left(  \mu^{2}+\left\vert A^{s}\right\vert ^{2}\right)
^{\frac{p}{2}}+\lambda\left(  \mu^{2}+\left\vert A^{s}\right\vert ^{2}%
+\beta^{2}\left\vert A^{a}\right\vert ^{2}\right)  ^{\frac{p}{2}}. \label{900}%
\end{equation}
It is clear that (\ref{5}) holds, while (\ref{6}) follows from (\ref{18}).

It remains to prove that, for some $\beta>0$, the function\ $F$ is
1-quasiconvex, i.e.,
\begin{equation}
\int_{Q}\left[  F\left(  A+\nabla\varphi\right)  -F\left(  A\right)  \right]
\,dx\geq0 \label{200}%
\end{equation}
for every $A\in\mathbb{M}^{n\times n}$ and for every $\varphi\in
C_{\operatorname*{per}}^{\infty}(Q;\mathbb{R}^{n})$.

Let $f_{\lambda}$ be the 2-quasiconvex function defined in Lemma
\ref{lemma 11}, and let
\[
g_{\beta}(A)=\hat{g}_{\beta}\left(  A^{s},A^{a}\right)  :=\left(  \mu
^{2}+\left|  A^{s}\right|  ^{2}+\beta^{2}\left|  A^{a}\right|  ^{2}\right)
^{\frac{p}{2}},
\]
so that
\[
F(A)=f_{\lambda}(A^{s})+\lambda g_{\beta}(A).
\]

Let us prove (\ref{200}). Fix a $Q$-periodic function $\varphi:\mathbb{R}%
^{n}\rightarrow\mathbb{R}^{n}$ of class $C^{\infty}$ and consider the periodic
Helmholtz decomposition
\[
\varphi=\nabla\phi+\psi
\]
given by Lemma \ref{lemma 6}. Then we have%
\begin{align}
\int_{Q}[F  &  \left(  A+\nabla\varphi\right)  -F\left(  A\right)
]\,dx\nonumber\\
=  &  \int_{Q}\left[  f_{\lambda}\left(  A^{s}+\nabla\varphi^{s}\right)
-f_{\lambda}\left(  A^{s}+\nabla\varphi^{s}-\nabla\psi^{s}\right)  \right]
\,dx\nonumber\\
&  +\int_{Q}\left[  f_{\lambda}\left(  A^{s}+\nabla^{2}\phi\right)
-f_{\lambda}\left(  A^{s}\right)  \right]  \,dx\label{170}\\
&  +\lambda\int_{Q}\left[  g_{\beta}\left(  A+\nabla\varphi\right)  -g_{\beta
}\left(  A\right)  \right]  \,dx\nonumber\\
=  &  \!:I_{1}+I_{2}+I_{3}.\nonumber
\end{align}
Since $\nabla f_{\lambda}\left(  A^{s}\right)  $ is a symmetric matrix, we
have $\nabla f_{\lambda}\left(  A^{s}\right)  \cdot\nabla\psi^{s}=\nabla
f_{\lambda}\left(  A^{s}\right)  \cdot\nabla\psi$, and therefore, by
periodicity,
\[
\int_{Q}\nabla f_{\lambda}\left(  A^{s}\right)  \cdot\nabla\psi^{s}\,dx=0.
\]
Hence
\[
I_{1}=-\int_{Q}\left[  f_{\lambda}\left(  A^{s}+\nabla\varphi^{s}-\nabla
\psi^{s}\right)  -f_{\lambda}\left(  A^{s}+\nabla\varphi^{s}\right)  +\nabla
f_{\lambda}\left(  A^{s}\right)  \cdot\nabla\psi^{s}\right]  \,dx.
\]
Since the function $g$ defined in (\ref{950}) clearly satisfies condition
(\ref{12}) , by Lemma \ref{Lemma 3} and (\ref{3}) it follows that (\ref{3})
still holds for the function $f_{\lambda}$ for a suitable constant $M>0$ in
place of $L$. We are now in position to apply Lemma \ref{Lemma 4} to obtain a
constant $\sigma=\sigma_{p,M}$ such that
\begin{align*}
I_{1}  &  \geq-\lambda\theta_{p}\int_{Q}\left(  \mu^{2}+\left\vert
A^{s}\right\vert ^{2}+\left\vert \nabla\varphi^{s}\right\vert ^{2}\right)
^{\frac{p-2}{2}}\left\vert \nabla\varphi^{s}\right\vert ^{2}\,dx\\
&  -\sigma\left(  \mu^{2}+\left\vert A^{s}\right\vert ^{2}\right)
^{\frac{p-2}{2}}\int_{Q}\left\vert \nabla\psi^{s}\right\vert ^{2}%
\,dx-\sigma\int_{Q}\left\vert \nabla\psi^{s}\right\vert ^{p}\,dx,
\end{align*}
and so, using Lemma \ref{lemma 1},\ we obtain%
\begin{align}
I_{1}  &  \geq-\lambda\theta_{p}\int_{Q}\left(  \mu^{2}+\left\vert
A^{s}\right\vert ^{2}+\left\vert \nabla\varphi^{s}\right\vert ^{2}\right)
^{\frac{p-2}{2}}\left\vert \nabla\varphi^{s}\right\vert ^{2}\,dx\label{201}\\
&  -\sigma\gamma_{n,2}\left(  \mu^{2}+\left\vert A^{s}\right\vert ^{2}\right)
^{\frac{p-2}{2}}\int_{Q}\left\vert \nabla\psi^{a}\right\vert ^{2}%
\,dx-\sigma\gamma_{n,p}\int_{Q}\left\vert \nabla\psi^{a}\right\vert
^{p}\,dx.\nonumber
\end{align}
On the other hand, the $2$-quasiconvexity of $f_{\lambda}$ yields
\begin{equation}
I_{2}\geq0. \label{202}%
\end{equation}
Since, by periodicity,
\[
\int_{Q}\nabla g_{\beta}(A)\cdot\nabla\varphi\,dx=0,
\]
by Lemma \ref{lemma 9} we have%
\begin{align}
I_{3}  &  =\lambda\int_{Q}\left[  g_{\beta}\left(  A+\nabla\varphi\right)
-g_{\beta}\left(  A\right)  -\nabla g_{\beta}(A)\cdot\nabla\varphi\right]
\,dx\nonumber\\
&  \geq\lambda\theta_{p}\int_{Q}\left(  \mu^{2}+\left\vert A^{s}\right\vert
^{2}+\left\vert \nabla\varphi^{s}\right\vert ^{2}\right)  ^{\frac{p-2}{2}%
}\left\vert \nabla\varphi^{s}\right\vert ^{2}\,dx\label{203}\\
&  +\frac{\lambda\theta_{p}\beta^{2}}{2}\left(  \mu^{2}+\left\vert
A^{s}\right\vert ^{2}\right)  ^{\frac{p-2}{2}}\int_{Q}\left\vert \nabla
\psi^{a}\right\vert ^{2}\,dx+\frac{\lambda\theta_{p}\beta^{p}}{2}\int
_{Q}\left\vert \nabla\psi^{a}\right\vert ^{p}\,dx.\nonumber
\end{align}
Choosing $\beta>0$ so that
\[
\frac{\lambda\theta_{p}\beta^{2}}{2}\geq\sigma\gamma_{n,2},\qquad\frac
{\lambda\theta_{p}\beta^{p}}{2}\geq\sigma\gamma_{n,p},
\]
by (\ref{201}), (\ref{202}), and (\ref{203}), we obtain
\[
I_{1}+I_{2}+I_{3}\geq0,
\]
which together with (\ref{170}) yields (\ref{200}).
\end{proof}

\begin{proof}
[Proof of Theorem \ref{theorem 1}]Since the function $t\mapsto f\left(
A+ta\otimes b+tb\otimes a\right)  $ is convex on $\mathbb{R}$ for every
$A\in\mathbb{M}_{\operatorname*{sym}}^{n\times n}$ and every $a$,
$b\in\mathbb{R}^{n}$ (see, e.g., [\cite{FM}]), from the growth condition (b)
it follows that there exists a constant $L>0$ depending only on $M$ and $p$
such that
\begin{equation}
\left\vert f\left(  A+B\right)  -f\left(  A\right)  \right\vert \leq L\left(
1+\left\vert A\right\vert ^{p-1}+\left\vert B\right\vert ^{p-1}\right)
\left\vert B\right\vert \label{114}%
\end{equation}
for every $A$, $B\in\mathbb{M}_{\operatorname*{sym}}^{n\times n}$. Given a
constant $\beta>0$, to be chosen at the end of the proof, let $G:\mathbb{M}%
^{n\times n}\rightarrow\mathbb{R}$ be the function defined by%
\begin{equation}
G(A):=f(A^{s})+\beta\left\vert A^{a}\right\vert ^{p}, \label{110}%
\end{equation}
and let $F$ be its $1$-quasiconvexification.

We want to prove that there exist two increasing sequences of positive numbers
$\left\{  \beta_{k}\right\}  $ and $\left\{  \lambda_{k}\right\}  $, depending
only on $k$, $p$, $\mu$, $\nu$, $M$, but not on the specific function
$f$,\ such that the corresponding functions $G_{k}$ satisfy%
\begin{equation}
\int_{Q}\left[  G_{k}\left(  A+\nabla\varphi\right)  -G_{k}\left(  A\right)
\right]  \,dx\geq-\frac{1}{k}\left\vert A^{s}\right\vert ^{p}-\lambda
_{k}\left\vert A^{a}\right\vert ^{p}-\frac{1}{k} \label{111}%
\end{equation}
for every $A\in\mathbb{M}^{n\times n}$ and for every $Q$-periodic function
$\varphi:\mathbb{R}^{n}\rightarrow\mathbb{R}^{n}$ of class $C^{\infty}$. This
will imply (see, e.g., \cite{D}) that for every $A\in\mathbb{M}^{n\times n}$
we have%
\begin{equation}
G_{k}(A)-\frac{1}{k}\left\vert A^{s}\right\vert ^{p}-\lambda_{k}\left\vert
A^{a}\right\vert ^{p}-\frac{1}{k}\leq F_{k}(A)\leq G_{k}(A) \label{112}%
\end{equation}
which yields (\ref{1}) and (\ref{7}) since $G_{k}(A)=f\left(  A\right)  $
whenever $A\in\mathbb{M}_{\operatorname*{sym}}^{n\times n}$.

Let us prove (\ref{111}). Fix a $Q$-periodic function $\varphi:\mathbb{R}%
^{n}\rightarrow\mathbb{R}^{n}$ of class $C^{\infty}$ and consider the periodic
Helmholtz decomposition
\[
\varphi=\nabla\phi+\psi
\]
given by Lemma \ref{lemma 6}. Then we have%
\begin{align}
\int_{Q}[G  &  \left(  A+\nabla\varphi\right)  -G\left(  A\right)
]\,dx\nonumber\\
=  &  \int_{Q}\left[  f\left(  A^{s}+\nabla^{2}\phi+\nabla\psi^{s}\right)
-f\left(  A^{s}+\nabla^{2}\phi\right)  \right]  \,dx\nonumber\\
&  +\int_{Q}\left[  f\left(  A^{s}+\nabla^{2}\phi\right)  -f\left(
A^{s}\right)  \right]  \,dx\label{113}\\
&  +\beta\int_{Q}\left[  \left|  A^{a}+\nabla\psi^{a}\right|  ^{p}-\left|
A^{a}\right|  ^{p}\right]  \,dx\nonumber\\
=  &  \!:I_{1}+I_{2}+I_{3}.\nonumber
\end{align}
By (\ref{114}) and by Cauchy's inequality, for every $\delta>0$ there exists a
constant $c_{\delta,p,L}>0$ such that
\begin{align*}
I_{1}  &  \geq-L\int_{Q}\left(  1+\left|  A^{s}+\nabla^{2}\phi\right|
^{p-1}+\left|  \nabla\psi^{s}\right|  ^{p-1}\right)  \left|  \nabla\psi
^{s}\right|  \,dx\\
&  \geq-\delta-\delta\left|  A^{s}\right|  ^{p}-\delta\int_{Q}\left|
\nabla^{2}\phi\right|  ^{p}\,dx-c_{\delta,p,L}\int_{Q}\left|  \nabla\psi
^{s}\right|  ^{p}\,dx.
\end{align*}
Hence, using Lemma \ref{lemma 1} we obtain
\[
I_{1}\geq-\delta-\delta\left|  A^{s}\right|  ^{p}-\delta\int_{Q}\left|
\nabla^{2}\phi\right|  ^{p}\,dx-c_{\delta,p,L}\gamma_{n,p}\int_{Q}\left|
\nabla\psi^{a}\right|  ^{p}\,dx.
\]
\ 

If $p\geq2$ then we have
\begin{align*}
I_{1}  &  \geq-\delta-\delta\left|  A^{s}\right|  ^{p} -\delta\int_{Q}\left(
\mu^{2}+\left|  A^{s}\right|  ^{2}+\left|  \nabla^{2}\phi\right|  ^{2}\right)
^{\frac{p-2}{2}}\left|  \nabla^{2}\phi\right|  ^{2}\,dx\\
&  -c_{\delta,p,L}\gamma_{n,p}\int_{Q}\left(  \left|  A^{a}\right|
^{2}+\left|  \nabla\psi^{a}\right|  ^{2}\right)  ^{\frac{p-2}{2}}\left|
\nabla\psi^{a}\right|  ^{2}\,dx.
\end{align*}

If $1<p\leq2$ then by Lemma \ref{lemma 12}, and by Lemma \ref{lemma10} with
$\mu=0$, we obtain for every $0<\varepsilon<1$
\begin{align*}
I_{1}  &  \geq-\delta\left(  1+\varepsilon\mu^{p}\right)  -\delta\left(
1+\varepsilon\right)  \left|  A^{s}\right|  ^{p}\\
&  -8 \delta\varepsilon^{\frac{p-2}{p}}\int_{Q}\left(  \mu^{2}+\left|
A^{s}\right|  ^{2}+\left|  \nabla^{2}\phi\right|  ^{2}\right)  ^{\frac{p-2}%
{2}}\left|  \nabla^{2}\phi\right|  ^{2}\,dx\\
&  -c_{\delta,p,L}\gamma_{n,p}\varepsilon^{\frac{p-2}{2}}\int_{Q}\left(
\left|  A^{a}\right|  ^{2}+\left|  \nabla\psi^{a}\right|  ^{2}\right)
^{\frac{p-2}{2}}\left|  \nabla\psi^{a}\right|  ^{2}\,dx-c_{\delta,p,L}%
\gamma_{n,p}\varepsilon^{\frac{p}{2}}\left|  A^{a}\right|  ^{p}.
\end{align*}
In both cases there exists a sequence of positive numbers $\left\{
\lambda_{k}\right\}  $, depending only on $p$, $\mu$, $\nu$, $M$, such that
for every $k$%
\begin{align}
I_{1}  &  \geq-\frac{1}{k}-\frac{1}{k}\left|  A^{s}\right|  ^{p}-\nu\int
_{Q}\left(  \mu^{2}+\left|  A^{s}\right|  ^{2}+\left|  \nabla^{2}\phi\right|
^{2}\right)  ^{\frac{p-2}{2}}\left|  \nabla^{2}\phi\right|  ^{2}%
\,dx\label{401}\\
&  -\lambda_{k}\int_{Q}\left(  \left|  A^{a}\right|  ^{2}+\left|  \nabla
\psi^{a}\right|  ^{2}\right)  ^{\frac{p-2}{2}}\left|  \nabla\psi^{a}\right|
^{2}\,dx-\lambda_{k}\left|  A^{a}\right|  ^{p}.\nonumber
\end{align}
The strict $2$-quasiconvexity of $f$ (condition (a)) yields%
\begin{equation}
I_{2}\geq\nu\int_{Q}\left(  \mu^{2}+\left|  A^{s}\right|  ^{2}+\left|
\nabla^{2}\phi\right|  ^{2}\right)  ^{\frac{p-2}{2}}\left|  \nabla^{2}%
\phi\right|  ^{2}\,dx. \label{402}%
\end{equation}
Since, by periodicity,
\[
\int_{Q}A^{a}\cdot\nabla\psi^{a}\,dx=0,
\]
by Lemma \ref{lemma 8} with $\mu=0$, we have%
\begin{align}
I_{3}  &  =\beta\int_{Q}\left[  \left|  A^{a}+\nabla\psi^{a}\right|
^{p}-\left|  A^{a}\right|  ^{p}-p\left|  A^{a}\right|  ^{p-2}A^{a}\cdot
\nabla\psi^{a}\right]  \,dx\label{403}\\
&  \geq\beta\theta_{p}\int_{Q}\left(  \left|  A^{a}\right|  ^{2}+\left|
\nabla\psi^{a}\right|  ^{2}\right)  ^{\frac{p-2}{2}}\left|  \nabla\psi
^{a}\right|  ^{2}\,dx.\nonumber
\end{align}
Choosing $\beta_{k}>0$ so that $\beta_{k}\theta_{p}\geq\lambda_{k}$, from
(\ref{401}), (\ref{402}), and (\ref{403}), we obtain
\[
I_{1}+I_{2}+I_{3}\geq-\frac{1}{k}-\frac{1}{k}\left|  A^{s}\right|
^{p}-\lambda_{k}\left|  A^{a}\right|  ^{p},
\]
which together with (\ref{113}) gives (\ref{111}).
\end{proof}

\begin{remark}
\label{remark1}\emph{It is clear from the proof of Theorem \ref{theorem 1}
that if }$f$\emph{ is nonnegative the we may take }$F_{k}$\emph{ to be also
nonnegative.}
\end{remark}

\section{Lower semicontinuity\label{section lower}}

The proof of Theorem \ref{theorem4} relies on the so-called Decomposition
Lemma (see \cite{FMP}).

\begin{lemma}
[Decomposition Lemma]\label{lemma decomposition}Let $\Omega$ be a bounded open
set in $\mathbb{R}^{n}$, let $p>1$, and let $\left\{  u_{k}\right\}  $ be a
sequence weakly converging to a function $u$ in $W^{1,p}\left(  \Omega
;\mathbb{R}^{n}\right)  $. Then there exists a subsequence (not relabeled) and
a sequence $\left\{  v_{k}\right\}  $ weakly converging to $u$ in
$W^{1,p}\left(  \Omega;\mathbb{R}^{n}\right)  $ such that $v_{k}=u$ in a
neighborhood of $\partial\Omega$, $\left\{  \left\vert \nabla v_{k}\right\vert
^{p}\right\}  $\ is equi-integrable, and $\mathcal{L}^{n}\left(  \left\{
u_{k}\neq v_{k}\right\}  \right)  \rightarrow0$.
\end{lemma}

The following simple lemma may be found in \cite{FLP}, however we include its
proof for the convenience of the reader.

\begin{lemma}
\label{lemma lower} Let $D\subset{\mathbb{R}}^{m}$ be an open set and let
\[
f:D\times\mathbb{M}^{d\times n}\rightarrow{\mathbb{R}}%
\]
be a lower semicontinuous function such that for every $v\in D$ the function
$f(v,\cdot)$ is continuous. Then for every $\bar{v}\in D$, $\varepsilon>0$,
and $L>0$ there exists $\delta=\delta(\bar{v},\varepsilon,L)\in(0,1)$ such
that
\[
f(\bar{v},A)\leq f(v,A)+\varepsilon
\]
for every $(v,A)\in D\times{\mathbb{M}}^{d\times n}$, with $|v-\bar{v}%
|\leq\delta$ and $|A|\leq L$.
\end{lemma}

\begin{proof}
Assume, for contradiction, that there exist $\bar{v}\in D$, $L>0$,
$\bar{\varepsilon}>0$, and a sequence
\[
\{(v_{k},A_{k})\}\subset D\times\overline{B_{d\times n}(0,L)},
\]
such that
\begin{equation}
\bar{\varepsilon}+f(v_{k},A_{k})<f(\bar{v},A_{k}) \label{not}%
\end{equation}
and $(v_{k},A_{k})\rightarrow(\bar{v},\bar{A})$ as $k\rightarrow\infty$, for
some $\bar{A}\in\overline{B_{d\times n}(0,L)}$. Since the function $f(\bar
{v},\cdot)$ is continuous and $f$ is lower semicontinuous, for any
$\varepsilon<\frac{1}{2}\bar{\varepsilon}$ there exists $\delta>0$ such that
\[
|f(\bar{v},A)-f(\bar{v},\bar{A})|\leq\varepsilon,\qquad f(\bar{v},\bar{A})\leq
f(v,A)+\varepsilon
\]
for all $(v,A)\in D\times\overline{B_{d\times n}(0,L)}$ with
\[
|v-\bar{v}|+|A-\bar{A}|\leq\delta.
\]
Thus for all $n$ sufficiently large, also by (\ref{not}), we have
\[
\bar{\varepsilon}+f(v_{k},A_{k})<f(\bar{v},A_{k})\leq f(\bar{v},\bar
{A})+\varepsilon\leq f(v_{k},A_{k})+2\varepsilon,
\]
which is a contradiction.
\end{proof}

Although the following lemma is well known to experts, its proof is not easy
to find in the literature and so we present it below for the reader's convenience.

\begin{lemma}
Let $D\subset{\mathbb{R}}^{m}$ be an open set and let
\[
f:D\times\mathbb{M}^{d\times n}\rightarrow{\mathbb{R}}%
\]
be a lower semicontinuous function which satisfies the following conditions:

\begin{enumerate}
\item[(a)] for every $v\in D$ the function $f(v,\cdot)$ is continuous in
$\mathbb{M}^{d\times n}$;

\item[(b)] there exist a
locally bounded
function $a:D\rightarrow\lbrack0,+\infty)$, a lower semicontinuous function
$b:D\rightarrow(0,+\infty)$, a locally bounded function $c:D\rightarrow
\lbrack0,+\infty)$, and a constant $p>1$ such that%
\[
b\left(  v\right)  |A|^{p}-c(v)\leq f(v,A)\leq a(v)(1+|A|^{p})\
\]
for every $(v,A)\in D\times\mathbb{M}^{d\times n}$.
\end{enumerate}

\noindent For every $v\in D$, let $\mathcal{Q}f(v,\cdot)$ be the
$1$-quasiconvexification of the function $f(v,\cdot)$. Then $\mathcal{Q}f$ is
lower semicontinuous on $D\times\mathbb{M}^{d\times n}$.
\end{lemma}

\begin{proof}
By conditions (a) and (b)\ for every $(v,A)\in D\times\mathbb{M}^{d\times n}$,
we have (see~\cite{D})
\begin{equation}
\mathcal{Q}f(v,A)=\inf\left\{  \int_{Q}f(v,A+\nabla\varphi(x))\,dx:\,\varphi
\in C_{c}^{1}(Q;\mathbb{R}^{d})\right\}  . \label{500}%
\end{equation}
By replacing $f(v,A)$ with $f(v,A)+\tilde{c}(v)$, where $\tilde{c}$\ is any
continuous function with $\tilde{c}\geq c$,\ we may assume without loss of
generality that $f\geq0$.

We begin by showing that for every fixed $A\in\mathbb{M}^{d\times n}$ the
function $\mathcal{Q}f(\cdot,A)$ is lower semicontinuous. Without loss of
generality we may assume that $A=0$. Let $\left\{  v_{k}\right\}  \subset D$
be a sequence converging to some $\bar{v}\in D$. If
\[
\liminf_{k\rightarrow\infty}\mathcal{Q}f(v_{k},0)=\infty,
\]
then there is nothing to prove. Thus, without loss of generality, we may
assume that
\[
\liminf_{k\rightarrow\infty}\mathcal{Q}f(v_{k},0)=\lim_{k\rightarrow\infty
}\mathcal{Q}f(v_{k},0)<\infty,
\]
and
\begin{equation}
C:=\sup_{k}\mathcal{Q}f(v_{k},0)<\infty. \label{501}%
\end{equation}
By (\ref{500}) for every fixed $0<\varepsilon<1$ and for every $k\in
\mathbb{N}$ there exists $\varphi_{k}\in C_{c}^{1}(Q;\mathbb{R}^{d})$\ such
that
\begin{equation}
\mathcal{Q}f(v_{k},0)+\varepsilon\geq\int_{Q}f(v_{k},\nabla\varphi
_{k}(x))\,dx. \label{502}%
\end{equation}
Hence, by condition (b) and (\ref{501}), we have
\begin{align*}
C+1  &  \geq b\left(  v_{k}\right)  \int_{Q}|\nabla\varphi_{k}(x)|^{p}%
\,dx-c(v_{k})\\
&  \geq b_{0}\int_{Q}|\nabla\varphi_{k}(x)|^{p}\,dx-c_{0},
\end{align*}
where $b_{0}:=\inf_{k}b(v_{k})>0$ and $c_{0}:=\sup_{k}c(v_{k})<\infty$, since
$b$ is lower semicontinuous and $c$ is locally bounded. By the Decomposition
Lemma (see Lemma \ref{lemma decomposition}), there exists a subsequence of
$\left\{  \varphi_{k}\right\}  $ (not relabeled) and a sequence $\left\{
w_{k}\right\}  $ weakly converging to some function $w$ in $W_{0}^{1,p}\left(
Q;\mathbb{R}^{d}\right)  $ such that $\left\{  \left\vert \nabla
w_{k}\right\vert ^{p}\right\}  $\ is equi-integrable, and $\mathcal{L}%
^{n}\left(  \left\{  \varphi_{k}\neq w_{k}\right\}  \right)  \rightarrow0$.
Hence we may find $L\geq1$ and $\bar{k}\in\mathbb{N}$ such that
\begin{equation}
\mathcal{L}^{n}\left(  \left\{  \varphi_{k}\neq w_{k}\right\}  \right)
+\int_{\left\{  \left\vert \nabla w_{k}\right\vert \geq L\right\}  }\left\vert
\nabla w_{k}\right\vert ^{p}\,dx\leq\varepsilon, \label{505}%
\end{equation}
for every $k\geq\bar{k}$. By Lemma \ref{lemma lower}\ there exists
$\delta=\delta(\bar{v},L,\varepsilon)\in(0,1)$ such that
\[
f(\bar{v},B)\leq f(v,B)+\varepsilon
\]
for all $(v,B)\in D\times\mathbb{M}^{d\times n}$, with $|v-\bar{v}|\leq\delta$
and all $|B|\leq L$. Therefore
\begin{align}
\mathcal{Q}f(v_{k},0)+\varepsilon &  \geq\int_{\left\{  \varphi_{k}%
=w_{k}\right\}  \cap\left\{  \left\vert \nabla w_{k}\right\vert <L\right\}
}f(v_{k},\nabla w_{k}(x))\,dx\nonumber\\
&  \geq\int_{\left\{  \varphi_{k}=w_{k}\right\}  \cap\left\{  \left\vert
\nabla w_{k}\right\vert <L\right\}  }f(\bar{v},\nabla w_{k}%
(x))\,dx-\varepsilon\label{507}\\
&  \geq\int_{\left\{  \varphi_{k}=w_{k}\right\}  }f(\bar{v},\nabla
w_{k}(x))\,dx-\varepsilon(2+a(\bar{v})),\nonumber
\end{align}
where in the last inequality we have used (\ref{505}) and condition (b). Since
$\left\{  \left\vert \nabla w_{k}\right\vert ^{p}\right\}  $\ is
equi-integrable and $\mathcal{L}^{n}\left(  \left\{  \varphi_{k}\neq
w_{k}\right\}  \right)  \rightarrow0$, by condition (b) we have
\[
\liminf_{k\rightarrow\infty}\int_{\left\{  \varphi_{k}=w_{k}\right\}  }%
f(\bar{v},\nabla w_{k}(x))\,dx=\liminf_{k\rightarrow\infty}\int_{Q}f(\bar
{v},\nabla w_{k}(x))\,dx,
\]
hence letting $k\rightarrow\infty$ in (\ref{507}) yields
\begin{align*}
\liminf_{k\rightarrow\infty}\mathcal{Q}f(v_{k},0)+\varepsilon &  \geq
\liminf_{k\rightarrow\infty}\int_{Q}f(\bar{v},\nabla w_{k}(x))\,dx-\varepsilon
(2+a(\bar{v}))\\
&  \geq\mathcal{Q}f(\bar{v},0)-\varepsilon(2+a(\bar{v})).
\end{align*}
It is now sufficient to let $\varepsilon\rightarrow0^{+}$ to conclude that
$v\mapsto\mathcal{Q}f(v,A)$ is lower semicontinuous.

Finally, we observe that the continuity of $A\mapsto\mathcal{Q}f(v,A)$ is an
immediate consequence of the quasiconvexity of $\mathcal{Q}f(v,\cdot)$. Since
the coefficient $a(v)$ in (b) is locally bounded, the functions $A\mapsto
\mathcal{Q}f(v,A)$ have the same local modulus of continuity when $v$ varies
in a compact subset of $D$. Since we have seen that $v\mapsto\mathcal{Q}%
f(v,A)$ is lower semicontinuous on $D$, this implies that $\mathcal{Q}f$ is
lower semicontinuous on $D\times\mathbb{M}^{d\times n}$.
\end{proof}

\begin{proof}
[Proof of Theorem \ref{theorem4}]Fix $u\in W^{1,1}(\Omega)$ and let
$\{u_{j}\}\subset SBH(\Omega)$ be any sequence converging to $u$ in
$W^{1,1}(\Omega)$ and satisfying (\ref{700}). Without loss of generality we
may assume that
\[
\liminf_{j\rightarrow\infty}\int_{\Omega}f(x,u_{j},\nabla u_{j},\nabla
^{2}u_{j})\,dx=\lim_{j\rightarrow\infty}\int_{\Omega}f(x,u_{j},\nabla
u_{j},\nabla^{2}u_{j})\,dx.
\]
Fix $\varepsilon\in(0,1)$ and let
\[
f_{\varepsilon}(x,u,\xi,A):=f(x,u,\xi,A)+\varepsilon\left\vert A\right\vert
^{p}.
\]
The function $f_{\varepsilon}$ satisfies all the conditions of Theorem
\ref{theorem 1}, hence there exists an increasing sequence $\left\{
F_{k,\varepsilon}\right\}  $ of functions $F_{k,\varepsilon}:\Omega
\times{\mathbb{R}}\times{\mathbb{R}}^{n}\times\mathbb{M}^{n\times
n}\rightarrow\lbrack0,+\infty)$ such that for $\mathcal{L}^{n}$ a.e.
$x\in\Omega$ and for every $\left(  u,\xi\right)  \in{\mathbb{R}}%
\times{\mathbb{R}}^{n}$ the function $F_{k,\varepsilon}\left(  x,u,\xi
,\cdot\right)  $ is $1$-quasiconvex,
\begin{align}
&  \lim_{k\rightarrow\infty}F_{k,\varepsilon}\left(  x,u,\xi,A\right)
=f_{\varepsilon}\left(  x,u,\xi,A\right)  \qquad\forall A\in\mathbb{M}%
_{\operatorname*{sym}}^{n\times n},\label{600}\\
0  &  \leq F_{k,\varepsilon}\left(  x,u,\xi,A\right)  \leq c_{k}%
(x,u,\xi)\left(  1+\left\vert A\right\vert ^{p}\right)  \qquad\forall
A\in\mathbb{M}^{n\times n}. \label{601}%
\end{align}
Let $C:=\sup_{j}\left\Vert \nabla^{2}u_{j}\right\Vert _{L^{p}}^{p}$. For every
fixed $\varepsilon\in(0,1)$ and $k\in\mathbb{N}$, we have
\begin{align}
\lim_{j\rightarrow\infty}\int_{\Omega}f(x,u_{j},\nabla u_{j},\nabla^{2}%
u_{j})\,dx  &  \geq\liminf_{j\rightarrow\infty}\int_{\Omega}f_{\varepsilon
}(x,u_{j},\nabla u_{j},\nabla^{2}u_{j})\,dx-\varepsilon C\label{903}\\
&  \geq\liminf_{j\rightarrow\infty}\int_{\Omega}F_{k,\varepsilon}%
(x,u_{j},\nabla u_{j},\nabla^{2}u_{j})\,dx-\varepsilon C,\nonumber
\end{align}
where we used the fact that $f_{\varepsilon}\geq F_{k,\varepsilon}$. We note
that, in view of the construction in the proof of Theorem \ref{theorem 1}, the
function $F_{k,\varepsilon}$ is defined as $1$-quasiconvexification of
$f_{\varepsilon}(x,u,\xi,A^{s})+\beta\left\vert A^{a}\right\vert ^{p}$, and
so, by the previous lemma, we have that $F_{k,\varepsilon}$ is a normal
integrand. Define
\[
G_{k,\varepsilon}:\Omega\times\left(  \mathbb{R}\times\mathbb{R}^{n}\right)
\times\left(  \mathbb{R}^{n}\times\mathbb{M}^{n\times n}\right)
\rightarrow\left[  0,\infty\right)
\]
as%
\[
G_{k,\varepsilon}\left(  x,\left(  w_{1},w_{2}\right)  ,\left(  \xi_{1}%
,\xi_{2}\right)  \right)  :=F_{k,\varepsilon}\left(  x,w_{1},w_{2},\xi
_{2}\right)  .
\]
Indeed, with $w_{j}:=(u_{j},\nabla u_{j})$ \ then
\[
\int_{\Omega}F_{k,\varepsilon}(x,u_{j},\nabla u_{j},\nabla^{2}u_{j})\,dx
\]
reduces to
\[
\int_{\Omega}G_{k,\varepsilon}(x,w_{j},\nabla w_{j})\,dx.
\]
It is clear that $\left\{  w_{j}\right\}  \subset SBV(\Omega,\mathbb{R\times
R}^{n})$ and that
\[
\sup_{j}\left(  \left\Vert \nabla w_{j}\right\Vert _{L^{p}}+\int_{S(w_{j}%
)}\theta(\left\vert \left[  w_{j}\right]  \right\vert )\,d\mathcal{H}%
^{n-1}\right)  <\infty.
\]
Hence we may apply Theorem 1.2 in \cite{K} to obtain
\begin{align*}
\lim_{j\rightarrow\infty}  &  \int_{\Omega}f(x,u_{j},\nabla u_{j},\nabla
^{2}u_{j})\,dx\geq\liminf_{j\rightarrow\infty}\int_{\Omega}G_{k,\varepsilon
}(x,w_{j},\nabla w_{j})\,dx-\varepsilon C\\
&  \geq\int_{\Omega}G_{k,\varepsilon}(x,w,\nabla w)\,dx-\varepsilon
C=\int_{\Omega}F_{k,\varepsilon}(x,u,\nabla u,\nabla^{2}u)\,dx-\varepsilon C,
\end{align*}
where we have used (\ref{903}), and $w:=(u,\nabla u)$.

By Lebesgue's Monotone Convergence Theorem, letting $k\rightarrow\infty$ in
the previous inequality and using (\ref{600}) gives
\[
\lim_{j\rightarrow\infty}\int_{\Omega}f(x,u_{j},\nabla u_{j},\nabla^{2}%
u_{j})\,dx\geq\int_{\Omega}f_{\varepsilon}(x,u,\nabla u,\nabla^{2}%
u)\,dx-\varepsilon C.
\]
It now suffices to let $\varepsilon\rightarrow0^{+}$.
\end{proof}

\section*{Acknowledgments}

The research of I. Fonseca was partially supported by the National Science
Foundation under Grant No. DMS--0103799. The work of Gianni Dal Maso is part
of the Project \textquotedblleft Calculus of Variations\textquotedblright%
\ 2002, supported by the Italian Ministry of Education, University, and Research.

This work was undertaken when G. Dal Maso visited the Center for Nonlinear
Analysis (NSF Grant No. DMS--9803791), Carnegie Mellon University, Pittsburgh,
PA, USA. The authors thank the Center for Nonlinear Analysis for its support
during the preparation of this paper.

I. Fonseca was engaged in several stimulating discussions on the subject of
this paper with B. Dacorogna during her visit to EPFL (Lausanne, Switzerland)
in Fall 2001. G. Leoni wishes to thank J. Kristensen for pointing out
reference \cite{MS}.

The authors wish to thank D. Kinderlehrer for many interesting conversations
on the subject of this paper.

\bigskip

\bigskip

\noindent S.I.S.S.A,

\noindent Trieste, Italy

\noindent dalmaso@sissa.it

\bigskip

\noindent Department of Mathematical Sciences,

\noindent Carnegie Mellon University,

\noindent Pittsburgh, PA 15213, USA

\noindent fonseca@andrew.cmu.edu

\noindent giovanni@andrew.cmu.edu

\noindent morini@andrew.cmu.edu

\end{document}